\documentclass[aos,MSNbibl,dvips]{arximspdf}
\usepackage{mathbh}
\usepackage{graphicx}
\doi{10.1214/13-AOS1104} 
\volume{41}
\issue{2}
\pubyear{2013}
\firstpage{923}
\lastpage{956}

\makeatletter
\newcommand{\rrVert}{\Vert}
\newcommand{\rrvert}{\vert}
\newcommand{\llVert}{\Vert}
\newcommand{\llvert}{\vert}
\def\cal{\mathcal}
\newcommand{\eqref}[1]{(\ref{#1})}
\def\ind{\mathbh{1}}
\def\argmin{\mathop{\operatorname{argmin}}}
\def\supp{\mathop{\operatorname{supp}}}

\newproclaim{ass}{Assumption}
\newtheorem{coro}{Corollary}[section]
\newtheorem{lem}{Lemma}[section]
\newtheorem{prop}{Proposition}[section]
\newtheorem{theo}{Theorem}[section]

\newcommand{\bY}{\mathbf{Y}}
\newcommand{\CC}{\mathbb{C}}
\newcommand{\E}{\mathbb{E}}
\newcommand{\EE}{\mathbb{E}}
\newcommand{\FF}{\cal{F}}
\newcommand{\GG}{\mathcal{G}}
\newcommand{\grad}{\nabla}
\newcommand{\MM}{\mathcal{M}}
\newcommand{\R}{\cal{R}}
\newcommand{\RR}{\mathbb{R}}
\newcommand{\Sc}{\mathbb{S}}
\newcommand{\ZZ}{\mathbb{Z}}

\renewcommand{\P}{\mathbb{P}}
\renewcommand{\Pr}{g}
\makeatother

\begin{document}
\begin{frontmatter}

\title{Minimax properties of Fr\'echet means of discretely sampled curves}
\runtitle{Minimax properties of Fr\'echet means}

\begin{aug}
\author[A]{\fnms{J\'er\'emie} \snm{Bigot}\ead[label=e1]{Jeremie.Bigot@isae.fr}\thanksref{t1}}
\and
\author[B]{\fnms{Xavier} \snm{Gendre}\corref{}\ead[label=e2]{Xavier.Gendre@math.univ-toulouse.fr}}
\thankstext{t1}{Supported by the French Agence Nationale de la
Recherche (ANR) under reference ANR-JCJC-SIMI1 DEMOS.}
\affiliation{DMIA--ISAE and Institut de Math\'ematiques de Toulouse}
\address[A]{DMIA--ISAE\\
10 avenue \'Edouard Belin\\
BP 54032\\
31055 Toulouse Cedex 4\\
France\\
\printead{e1}}
\address[B]{Institut de Math\'ematiques de Toulouse\\
IMT---Universit\'e Paul Sabatier\\
118 route de Narbonne\\
31062 Toulouse Cedex 4\\
France\\
\printead{e2}}

\runauthor{J. Bigot and X. Gendre}
\end{aug}

\received{\smonth{10} \syear{2012}}
\revised{\smonth{2} \syear{2013}}

%
\begin{abstract}
We study the problem of estimating a mean pattern from a set of similar curves
in the setting where the variability in the data is due to random geometric
deformations and additive noise. We propose an estimator based on the
notion of
Fr\'echet mean that is a generalization of the standard notion of
averaging to
non-Euclidean spaces. We derive a minimax rate for this estimation
problem, and
we show that our estimator achieves this optimal rate under the asymptotics
where both the number of curves and the number of sampling points go to
infinity.
\end{abstract}

%

\end{frontmatter}

\section{Introduction}

\subsection{Fr\'echet means}

The Fr\'echet mean~\cite{fre} is an extension of the usual Euclidean
mean to nonlinear spaces endowed with non-Euclidean metrics. If $\bY
_{1},\ldots,\bY_{J}$ denote i.i.d. random variables with values in a
metric space $\MM$ with metric $d_{\MM}$, then the empirical Fr\'echet
mean $\overline{\bY}_{\MM}$ of the sample $\bY_{1},\ldots,\bY_{J}$ is
defined as a minimizer (not necessarily unique) of
\[
\overline{\bY}_{\MM} \in\argmin_{y \in\MM} \frac{1}{J} \sum
_{j=1}^{J} d_{\MM}^{2}(y,
\bY_{j}).
\]
For random variables belonging to a nonlinear manifold, a well-known
example is the computation of the mean of a set of planar shapes in
Kendall's shape space~\cite{kendall} that leads to the Procrustean
means studied in~\cite{MR1108330}. A detailed study of some properties
of the Fr\'echet mean in finite dimensional Riemannian manifolds (such
as consistency and uniqueness) has been performed in \cite
{Afsari,batach1,batach2,MR2816349}. However, there is not so much work
on the properties of the Fr\'echet mean in infinite dimensional and
non-Euclidean spaces of curves or images. In this paper, we are
concerned with the nonparametric estimation of a mean pattern
(belonging to a nonlinear space) from a set of similar curves in the
setting where the variability in the data is due to random geometric
deformations and additive noise.

More precisely, let us consider noisy realizations of $J$ curves $f_1,
\ldots, f_J \dvtx\break  [0,1] \to\RR$ sampled at $n$ equispaced points $t_{\ell}
= \frac{\ell}{n}, \ell=1,\ldots,n$,
%
\begin{equation}
\label{eqmodel} Y_{\ell,j} = f_{j}(t_{\ell}) +
\varepsilon_{\ell,j},\qquad \ell=1,\ldots,n \mbox{ and } j=1,\ldots,J,
\end{equation}
where the $\varepsilon_{\ell,j}$'s are independent and identically
distributed (i.i.d.) Gaussian variables with zero expectation and known
variance $\sigma^{2} > 0$. In many applications, the observed curves
have a similar structure that may lead to the assumption that the
$f_{j}$'s are random elements varying around the same mean pattern $f \dvtx [0,1] \to\RR$
(also called reference template). However, due to
additive noise and geometric variability in the data, this mean pattern
is typically unknown and has to be estimated. In this setting, a widely
used approach is Grenander's pattern theory~\cite{Gre,gremil,MR2176922,TYShape}
that models geometric variability by the
action of a Lie group on an infinite dimensional space of curves (or images).

When the curves $f_{j}$ in \eqref{eqmodel} exhibit a large source of
geometric variation in time, this may significantly complicates the
construction of a consistent estimator of a mean pattern. In what
follows, we consider the simple model of randomly shifted curves that
is commonly used in many applied areas such as neurosciences \cite
{MR2816476} or biology~\cite{MR1841413}. In such a framework, we have
%
\begin{equation}
\label{eqmodelshift} f_{j}(t) = f\bigl(t - \theta_{j}^{\ast}
\bigr)\qquad \mbox{for all } t \in[0,1] \mbox { and } j=1,\ldots,J,
\end{equation}
where $f \dvtx [0,1] \to\RR$ is an unknown curve that can be extended
outside $[0,1]$ by $1$-periodicity. In a similar way, we could consider
a function $f$ defined on the circle $\RR/\ZZ$. The shifts $\theta
_{j}^{\ast}$'s are supposed to be i.i.d. real random variables
(independent of the $\varepsilon_{\ell,j}$'s) that are sampled from an
unknown distribution $\mathbb{P}_{\ast}$ on $\RR$. In model~(\ref
{eqmodelshift}), the shifts $\theta_{j}^{\ast}$ represent a source of
geometric variability in time.

In functional data analysis, the problem of estimating a mean pattern
from a set of curves that differ by a time transformation is usually
referred to as the curve registration problem; see, for example, \cite
{ramsil}. Registering functional data has received a lot of attention
in the literature over the two last decades; see, for example, \cite
{big,kneipgas,ramsil,tangmull,wanggas} and references therein.
Nevertheless, in these papers, constructing consistent estimators of
the mean pattern $f$ as the number $J$ of curves tends to infinity is
generally not considered. Self-modeling regression methods proposed in
\cite{kg} are semiparametric models for curve registration that are
similar to the shifted curves model, where each observed curve is a
parametric transformation of an unknown mean pattern. Constructing a
consistent estimator of the mean pattern in such models has been
investigated in~\cite{kg} in an asymptotic framework where both the
number $J$ of curves and the number $n$ of design points grow toward
infinity. However, deriving optimal estimators in the minimax sense has
not been considered in~\cite{kg}. Moreover, a novel contribution of
this paper is to make a connection between the curve registration
problem and the notion of Fr\'echet mean in non-Euclidean spaces which
has not been investigated so far.

\subsection{Model and objectives}

The main goal of this paper is to construct nonparametric estimators of
the mean pattern $f$ from the data
%
\begin{equation}
\label{eqmodelshiftedcurve} Y_{\ell,j} = f\bigl(t_{\ell} -
\theta_{j}^{\ast}\bigr) + \varepsilon_{\ell,j},\qquad \ell=1,
\ldots,n \mbox{ and } j=1,\ldots,J,
\end{equation}
in the setting where both the number $J$ of curves and the number $n$
of design points are allowed to vary and to tend toward infinity.

In the sequel to this paper, it will be assumed that the random shifts
are sampled from an unknown density $g$ with respect to the Lebesgue
measure $d \theta$ [namely $d \mathbb{P}_{\ast}(\theta) = g(\theta) \,d
\theta$]. Note that since $f$ is assumed to be $1$-periodic, one may
restrict to the case where the density $g$ has a compact support
included in the interval $ [-\frac{1}{2},\frac{1}{2} ]$. Under
assumption \eqref{eqmodelshift}, the (standard) Euclidean mean $\bar
{Y}_{\ell} = \frac{1}{J} \sum_{j=1}^J Y_{\ell,j}$ of the data is
generally not a consistent estimator of the mean pattern $f$ at $t =
t_{\ell}$. Indeed, the law of large numbers implies that
\[
\lim_{J\to\infty} \bar{Y}^{\ell} = \lim_{J\to\infty}
\frac{1}{J} \sum_{j=1}^J f
\bigl(t_\ell- \theta_j^{\ast}\bigr) = \int
_{\RR} f(t_{\ell} - \theta )g(\theta)\,d\theta \qquad\mbox{a.s.} \]
Thus, under mild assumptions on $f$ and $g$, we have
\[
\lim_{J\to\infty} \bar{Y}^{\ell} =f\star g(t_{\ell})
\neq f(t_{\ell}) \qquad\mbox{a.s.},
\]
where $f\star g$ is the convolution product between $f$ and $g$.

To build a consistent estimator of $f$ in model \eqref
{eqmodelshiftedcurve}, we propose to use a notion of empirical Fr\'
echet mean in an infinite dimensional space. Recently, some properties
of Fr\'echet means in randomly shifted curves models have been
investigated in~\cite{MR2676894} and~\cite{BC11}. However, studying the
rate of convergence and the minimax properties of such estimators in
the double asymptotic setting $\min(n,J) \to+ \infty$ has not been
considered so far.

Note that model \eqref{eqmodelshiftedcurve} is clearly not
identifiable, as for any $\tilde{\theta} \in [-\frac{1}{2},\frac
{1}{2} ]$, one may replace $f(\cdot)$ by $\tilde{f}(\cdot) =f(\cdot
-\tilde{\theta})$ and $\theta_{j}^{\ast}$ by $\tilde{\theta}_{j} =
\theta_{j}^{\ast} - \tilde{\theta}$ without changing model~\eqref
{eqmodelshiftedcurve}. Therefore, estimation of $f$ is only feasible
up to a time shift. Thus, we propose to consider the problem of
estimating its equivalence class $[f]$ (or orbit) under the action of
shifts. More precisely, let $L^{2}_{\mathrm{per}}([0,1])$ be the space of
squared integrable functions on $[0,1]$ that can be extended outside
$[0,1]$ by 1-periodicity. Let $\Sc_{1}$ be the one-dimensional torus.
We recall that any element $\tau= \tau(\theta) \in\Sc_{1}$ can be
identified with an element $\theta$ in the interval $ [-\frac
{1}{2},\frac{1}{2} ]$. For $f \in L^{2}_{\mathrm{per}}([0,1])$, we define
its equivalence class by the action of a time shift as
\[
[f]:= \bigl\{ f^{\tau}, \tau\in\Sc_{1} \bigr\},
\]
where for $\tau= \tau(\theta) \in\Sc_{1} $ (with $\theta\in
[-\frac{1}{2},\frac{1}{2} ]$), $f^{\tau}(t) = f(t - \theta)$ for
all $t \in[0,1]$. Let $f,h \in L^{2}_{\mathrm{per}}([0,1])$, and we define the
distance between $[f], [h] \in L^{2}_{\mathrm{per}}([0,1]) /\break \Sc_{1} $ as
%
\begin{equation}
d\bigl([f],[h]\bigr) = \inf_{ \theta\in [-{1}/{2},{1}/{2} ] } \biggl( \int
_{0}^{1} \bigl|f(t-\theta) - h(t)\bigr|^{2}\,dt
\biggr)^{1/2}. \label {defdistance}
\end{equation}
In the setting of Grenander's pattern theory, $(L^{2}_{\mathrm{per}}([0,1]) / \Sc
_{1},d)$ represents an infinite dimensional and nonlinear set of
curves, and $\Sc_{1}$ is a Lie group modeling geometric variability in
the data.

\subsection{Main contributions}

Let us assume that $\FF\subset L^{2}_{\mathrm{per}}([0,1])$ represents some
smoothness class of functions (e.g., a Sobolev ball). Suppose also that
the unknown density $g$ of the random shifts in \eqref{eqmodelshift}
belongs to some set $\GG$ of probability density functions on $
[-\frac{1}{2},\frac{1}{2} ]$. Let $\hat{f}_{n,J}$ be some
estimator of $f$ based on the random variables $Y_{\ell,j}$ given by
\eqref{eqmodelshiftedcurve} taking its values in
$L^{2}_{\mathrm{per}}([0,1])$. For some $f \in\FF$, the risk of the estimator
$\hat{f}_{n,J}$ is defined by
\[
\R_{\Pr}(\hat{f}_{n,J},f) = \EE^{\Pr} \bigl(
d^{2}\bigl([\hat{f}_{n,J}],[f]\bigr) \bigr),
\]
where the above expectation $\EE^{\Pr} $ is taken with respect to the
distribution of the $Y_{\ell,j}$'s in \eqref{eqmodelshiftedcurve} and
under the assumption that the shifts are i.i.d. random variables
sampled from the density $\Pr$. We propose to investigate the
optimality of an estimator by introducing the following minimax risk:
\[
\R_{n,J}(\FF,\GG) = \inf_{\hat{f}_{n,J}} \sup
_{\Pr\in\GG} \sup_{f \in
\FF} \R_{\Pr}(
\hat{f}_{n,J},f),
\]
where the above infimum is taken over the set of all possible
estimators in model~(\ref{eqmodelshiftedcurve}).

For $f \in L^{2}_{\mathrm{per}}([0,1])$, let us denote its Fourier coefficients by
\[
c_{k} = \int_{0}^{1} f(t)
e^{- i 2 \pi k t} \,dt,\qquad k \in\ZZ.
\]
Suppose that $\FF= \tilde{W}_{s}(A,c_*)$ is the following bounded set
of nonconstant functions with degree of smoothness $s > 1/2$:
\[
\tilde{W}_{s}(A,c_*) = \biggl\{ f \in L^{2}_{\mathrm{per}}
\bigl([0,1]\bigr); \sum_{k \in
\ZZ} \bigl(1+|k|^{2s }
\bigr) |c_{k}|^2 \leq A^2\mbox{ with } \vert
c_{1}\vert\geq c_* \biggr\}
\]
for some positive reals $A$ and $c_*$. The introduction of the above
set is motivated by the definition of Sobolev balls. The additional
assumption $\vert c_1\vert\geq c_*$ is needed to ensure identifiability
of $[f]$ in model \eqref{eqmodelshiftedcurve} with respect to the
distance \eqref{defdistance}.

Moreover, let $\GG^{\kappa}$ be a set of probability densities having a
compact support of size smaller than $\kappa$ with $0 < \kappa< 1/8$
defined as
\[
\GG^{\kappa} = \biggl\{g \dvtx \biggl[-\frac{1}{2},\frac{1}{2}
\biggr]\to\RR ^{+}; \int_{-1/2}^{1/2}g(
\theta) \,d\theta= 1 \mbox{ and } \supp(g) \subseteq[-\kappa/2,\kappa/2] \biggr\}.
\]
Suppose also that the following condition holds:
\[
J \asymp n^{\alpha}\qquad \mbox{ for some } 0 < \alpha< 1/6,
\]
where the notation $J \asymp n^{\alpha}$ means that there exist two
positive constants $c_{2} \geq c_{1} > 0$ such that $ c_{1} n^{\alpha}
\leq J \leq c_{2} n^{\alpha}$ for any choices of $J$ and $n$.

Then, under such assumptions, the main contribution of the paper is to
show that one can construct an estimator $\hat{f}_{n,J} $ based on a
smoothed Fr\'echet mean of discretely sampled curves that satisfies
\[
\limsup_{ \min(n,J) \to+ \infty} r_{n,J}^{-1} \sup
_{\Pr\in\GG^{\kappa
}} \sup_{f \in\tilde{W}_s(A,c_*)} \R_{\Pr}(
\hat{f}_{n,J},f) \leq C_{0},
\]
where $C_{0} > 0$ is a constant that only depends on $A$, $s$, $\kappa
$, $c_*$ and $\sigma^2$. The rate of convergence $r_{n,J}$ is given by
\[
r_{n,J} = n^{-1} + (nJ)^{-{2s}/{(2s +1)}}.
\]
The two terms in the rate $r_{n,J}$ have different interpretations. The
second term $(nJ)^{-{2s}/{(2s +1)}}$ is the usual nonparametric rate
for estimating the function $f$ (over a Sobolev ball) in model \eqref
{eqmodelshiftedcurve} that we would obtain if the true shifts $\theta
_{1}^{\ast},\ldots,\theta_{J}^{\ast}$ were known. Moreover, under some
additional assumptions, we will show that this rate is optimal in the
minimax sense and that our estimator achieves it.

The first term $n^{-1}$ in the rate $r_{n,J}$ can be interpreted as
follows. As shown later in the paper, the computation of a Fr\'echet
mean of curves is a two-step procedure. It consists of building
estimators $\hat{\theta}_{j}^{0}$ of the unknown shifts and then
aligning the observed curves. For $\theta= (\theta_{1},\ldots,\theta
_{J}) \in\RR^{J}$, let us define the Euclidean norm $\| \theta\| =
( \sum_{j=1}^{J} |\theta_{j}|^{2} )^{1/2}$. One of the
contributions of this paper is to show that estimation of the vector
\[
\theta^{0} = \bigl(\theta_{1}^{\ast} - \bar{
\theta}_J, \ldots, \theta _{J}^{\ast}- \bar{
\theta}_J \bigr)' \in\RR^{J} \qquad\mbox{where }
\bar {\theta}_J=\frac{1}{J}\sum_{j=1}^{J}
\theta_j^{\ast},
\]
is feasible at the rate $n^{-1}$ for the normalized quadratic risk
$\frac{1}{J}\EE\| \hat{\theta}^{0} - \theta^{0} \|^{2} $ and that this
allows us to build a consistent Fr\'echet mean. If the number $J$ of
curves was fixed, $n ^{-1}$ would correspond to the usual
semi-parametric rate for estimating the shifts in model \eqref
{eqmodelshiftedcurve} in the setting where the $\theta^{\ast}_{j}$
are nonrandom parameters; see~\cite{BLV09,MR2369028,MR2662362} for
further details. Here, this rate of convergence has been obtained in
the double asymptotic setting $\min(n,J) \to+ \infty$. This setting
significantly complicates the estimation of the vector $\theta^{0} \in
\RR^{J}$ since its dimension $J$ is increasing with the sample size
$nJ$. Hence, in the case of $\min(n,J) \to+ \infty$, estimating the
shifts at the rate $n^{-1}$ is not a standard semi-parametric problem,
and we have to impose the constraint $J \asymp n^{\alpha}$ (with $0 <
\alpha< 1/6$) to obtain this result. The term $n^{-1}$ in the rate
$r_{n,J}$ is thus the price to pay for not knowing the random shifts in~\eqref{eqmodelshiftedcurve} that need to be estimated to compute a
Fr\'echet mean.

\subsection{Organization of the paper}

In Section~\ref{secsmoothedFrechet}, we introduce a notion of smoothed
Fr\'echet means of curves. We also discuss the connection between this
approach and the well-known problems of curve registration and image
warping. In Section~\ref{equpperbound}, we discuss the rate of
convergence of the estimators of the shifts. We also build a Fr\'echet
mean using model selection techniques, and we derive an upper bound on
its rate of convergence. In Section~\ref{eqlowerbound}, we derive a
lower bound on the minimax risk, and we give some sufficient conditions
to obtain a smoothed Fr\'echet mean converging at an optimal rate in
the minimax sense. In Section~\ref{secdiscussion}, we discuss the main
results of the paper and their connections with the nonparametric
literature on deformable models. Some numerical experiments on
simulated data are presented in Section~\ref{secsimus}. The proofs of
the main results are gathered in a technical \hyperref[app]{Appendix}.

\section{Smoothed Fr\'echet means of curves} \label{secsmoothedFrechet}

Let $f_{1},\ldots,f_{J} $ be a set of functions in
$L^{2}_{\mathrm{per}}([0,1])$. We define the Fr\'echet mean $[\bar{f}]$ of
$[f_{1}],\ldots,[f_{J}]$ as
\[
[\bar{f}] \in\argmin_{ [f] \in L^{2}_{\mathrm{per}}([0,1])/\Sc_{1}} \frac{1}{J} \sum
_{j=1}^{J} d^{2}\bigl([f],[f_{j}]
\bigr).
\]
It can be easily checked that a representation $\bar{f} \in
L^{2}_{\mathrm{per}}([0,1])$ of the class $[\bar{f}]$ is given by the following
two-step procedure:
\begin{longlist}[(1)]
\item[(1)] Computation of shifts to align the curves
%
\begin{eqnarray}
\label{eqstep1} \qquad&&(\tilde{\theta}_{1},\ldots,\tilde{
\theta}_{J})
\nonumber
\\[-8pt]
\\[-8pt]
\nonumber
&&\qquad\in\argmin_{ (\theta
_{1},\ldots,\theta_{J}) \in [-{1}/{2},{1}/{2} ]^{J} } \frac{1}{J} \sum
_{j=1}^{J} \int_{0}^{1}
\Biggl\llvert f_{j}(t+\theta_{j}) - \frac{1}{J}
\sum_{j'=1}^{J} f_{j'}(t+
\theta_{j'}) \Biggr\rrvert ^{2}\,dt.
\end{eqnarray}
\item[(2)] Averaging after an alignment step:
$
\bar{f}(t) = \frac{1}{J} \sum_{j=1}^{J} f_{j}(t+\tilde{\theta}_{j}),
t \in[0,1].
$
\end{longlist}

Let us now explain how the above two-step procedure can be used to
define an estimator of $f$ in model \eqref{eqmodelshiftedcurve}. Let
\[
\bigl\{ \phi_{k}(t) = e^{i 2 \pi k t}, t \in[0,1] \bigr
\}_{ k \in\ZZ}
\]
be the standard Fourier basis. For legibility, we assume that $n = 2 N
\geq4$ is even, and we split the data into two samples as follows:
\[
Y_{q,j}^{(0)} = Y_{2q,j}\quad \mbox{and}\quad
Y_{q,j}^{(1)} = Y_{2q-1,j}, \qquad q =1,\ldots,N
\]
for $j=1,\ldots,J$, and
\[
t_{q}^{(0)} = t_{2q}\quad \mbox{and}\quad
t_{q}^{(1)} = t_{2q-1}, \qquad q =1,\ldots,N.
\]
For any $z \in\CC$, we denote by $\overline{z}$ its complex conjugate.
Then we define the following empirical Fourier coefficients, for any
$j\in\{1,\ldots,J\}$:
\begin{eqnarray*}
\hat{c}_{k,j}^{(0)} & = & \frac{1}{N} \sum
_{q = 1 }^{N } Y_{q,j}^{(0)}
\overline{\phi_{k}\bigl(t_{q}^{(0)}\bigr)} =
\bar{c}_{k,j}^{(0)} + \frac{1}{\sqrt {N}} z_{k,j}^{(0)},\qquad
- \frac{N}{2} \leq k < \frac{N}{2},
\\
\hat{c}_{k,j}^{(1)} & = & \frac{1}{N} \sum
_{q = 1}^{N} Y_{q,j}^{(1)}
\overline{\phi_{k}\bigl(t_{q}^{(1)}\bigr)} =
\bar{c}_{k,j}^{(1)} + \frac{1}{\sqrt {N}} z_{k,j}^{(1)},\qquad
- \frac{N}{2} \leq k < \frac{N}{2},
\end{eqnarray*}
where
\[
\bar{c}_{k,j}^{(p)}= \frac{1}{N} \sum
_{q = 1}^{N} f\bigl(t_{q}^{(p)} -
\theta _{j}^{\ast}\bigr) \overline{\phi_{k}
\bigl(t_{q}^{(p)}\bigr)}, \qquad p \in\{0,1\},
\]
and the $z_{k,j}^{(p)}$'s are i.i.d. complex Gaussian variables with
zero expectation and variance $\sigma^{2}$.

Then we define estimators of the unknown random shifts $\theta_{j}^{\ast
}$ as
%
\begin{equation}
(\hat{\theta}_{1},\ldots,\hat{\theta}_{J}) \in
\argmin_{ (\theta
_{1},\ldots,\theta_{J}) \in [-{1}/{2},{1}/{2} ]^{J} } M_{n}(\theta_{1},\ldots,
\theta_{J}), \label{eqcrittheta}
\end{equation}
where
%
\begin{equation}
M_{n}(\theta_{1},\ldots,\theta_{J}) =
\frac{1}{J} \sum_{j=1}^{J} \sum
_{|k| \leq k_{0}} \Biggl\llvert \hat{c}_{k,j}^{(0)}
e^{i 2 \pi k \theta_{j}} - \frac{1}{J} \sum_{j'=1}^{J}
\hat{c}_{k,j'}^{(0)} e^{i 2 \pi k \theta
_{j'}} \Biggr\rrvert
^{2},\label{defMn}
\end{equation}
with some positive integer $k_{0}$ that will be discussed later. The
smoothed Fr\'echet mean of $f$ is then defined as
\[
\hat{f}_{n,J}^{(m)}(t) = \sum_{ |k| \leq m }
\Biggl( \frac{1}{J} \sum_{j=1}^{J}
\hat{c}_{k,j}^{(1)} e^{i 2 \pi k \hat{\theta}_{j}} \Biggr)
\phi_{k}(t) = \frac{1}{J} \sum_{j=1}^{J}
\hat{f}^{(m)}_{j}(t + \hat {\theta}_{j}), \qquad t
\in[0,1],
\]
where the integer $m\in\{1,\ldots,N/2\}$ is a frequency cut-off
parameter that will be discussed later and $\hat{f}^{(m)}_{j}(t) = \sum_{ |k| \leq m } \hat{c}_{k,j}^{(1)} \phi_{k}(t)$. Note that the
estimators\vadjust{\goodbreak} $\hat{\theta}_{j}$ of the shifts have been computed using
only half of the data and that the curves $\hat{f}^{(m)}_{j}$ are
calculated using the other half of the data. By splitting the data in
such a way, the random variables $\hat{\theta}_{j}$ and $\hat
{f}^{(m)}_{j}$ are independent conditionally to $(\theta_{1}^{\ast
},\ldots,\theta_{J}^{\ast})$. The computation of the $\hat{\theta}_j$'s
can be performed by using a gradient descent algorithm to minimize the
criterion \eqref{defMn}; for further details, see~\cite{MR2676894}.

Note also that this two-step procedure does not require the use of a
reference template to compute estimators $\hat{\theta}_{1},\ldots,\hat
{\theta}_{J}$ of the random shifts. Indeed, one can interpret the term
$\frac{1}{J} \sum_{j'=1}^{J}\hat f_{j'}(t+\theta_{j'})$ in \eqref
{eqstep1} as a template that is automatically estimated. In
statistics, estimating a mean pattern from set of curves that differ by
a time transformation is usually referred to as the curve registration
problem. It has received a lot of attention over the last two decades;
see, for example,~\cite{big,MR1841413,MR2816476} and
references therein. Hence, there exists a connection between our
approach and the well-known problems of curve registration and its
generalization to higher dimensions (image warping); see, for example,
\cite{MR1858399}. However, studying the minimax properties of an
estimator of a mean pattern in curve registration models has not been
investigated so far.

\section{Upper bound on the risk} \label{equpperbound}

\subsection{Consistent estimation of the unknown shifts}

Note that, due to identifiability issues in model \eqref
{eqmodelshiftedcurve}, the minimization (\ref{eqcrittheta}) is not
well defined. Indeed, for any $(\hat{\theta}_{1},\ldots,\hat{\theta
}_{J})$ that minimizes \eqref{defMn}, one has that for any $\tilde
{\theta}$ such that $(\hat{\theta}_{1}+\tilde{\theta},\ldots,\hat{\theta
}_{J}+\tilde{\theta}) \in [-\frac{1}{2},\frac{1}{2} ]^{J} $,
this vector is also a minimizer of $M_{n}$. Choosing identifiability
conditions amounts to imposing constraints on the minimization of the criterion
%
\begin{equation}
M(\theta_{1},\ldots,\theta_{J}) = \frac{1}{J} \sum
_{j=1}^{J} \sum
_{|k|
\leq k_{0}} \Biggl\llvert c_{k} e^{i2\pi k (\theta_{j}-\theta_{j}^{\ast})} -
\frac{1}{J} \sum_{j'=1}^{J}
c_{k} e^{i2\pi k (\theta_{j'}-\theta
_{j'}^{\ast})} \Biggr\rrvert ^{2}, \label{defM}
\end{equation}
where $c_{k}$, $k \in\ZZ$, are the Fourier coefficients of the mean
pattern $f$. Criterion \eqref{defM} can be interpreted as a version
without noise criterion \eqref{defMn} when replacing $\hat
{c}_{k,j}^{(0)}$ by $c_{k} e^{-i 2 \pi k \theta_j^{\ast}}$. Obviously,
criterion \eqref{defM} admits a minimum at $\theta^{\ast} = (\theta
_{1}^{\ast},\ldots,\theta_{J}^{\ast})$ such that $M(\theta^{\ast}) =
0$. However, this minimizer over $ [-\frac{1}{2},\frac{1}{2}
]^{J}$ is clearly not unique. To impose uniqueness of some minimum of
$M$ over a restricted set, let us introduce the following
identifiability conditions:
%
\begin{ass}\label{hypP}
The distribution $\Pr$ of the random shifts has a compact support
included in $[-\kappa/2,\kappa/2]$ for some $0 < \kappa< 1/8$.
\end{ass}
%
\begin{ass}\label{hypf}
The mean pattern $f$ in model \eqref{eqmodelshiftedcurve} is such
that $c_{1} = \int_0^1 f(x) e^{-i2\pi x}\,dx \neq0$.\vadjust{\goodbreak}
\end{ass}

Assumption~\ref{hypP} means that the support of the density $g$ of the
random shifts should be sufficiently small. This implies that the
shifted curves $f(t-\theta_{j}^{\ast})$ are somehow concentrated around
the unknown mean pattern $f$. Such an assumption of concentration of
the data around a reference shape has been used in various papers to
prove the uniqueness and the consistency of Fr\'echet means for random
variables lying in a Riemannian manifold; see \cite
{Afsari,batach1,batach2,MR2816349}. Assumption~\ref{hypP} could
certainly be weakened by dealing with a basis other than the Fourier
polynomials. However, this is not the main point of this paper. Recent
studies in this direction have been made to derive asymptotic on the
Fr\'echet mean of a distribution on the circle without restriction on
its support; see~\cite{Hotz1374530,Charlier} or~\cite{Mc2010},
for instance. Assumption~\ref{hypf} is an identifiability condition to
avoid the case where the function $f$ is constant over $[0,1]$ which
would make impossible the estimation of the unobserved random shifts.

For $0 < \kappa< 1/8$, let us define the constrained set
\[
\Theta_{\kappa} = \Biggl\{ (\theta_{1},\ldots,
\theta_{J}) \in[-\kappa /2,\kappa/2]^{J}, \sum
_{j=1}^{J} \theta_{j} = 0 \Biggr\}.
\]
Let
\[
\theta_{j}^{0} = \theta_{j}^{\ast} -
\frac{1}{J} \sum_{j'=1}^{J} \theta
_{j'}^{\ast},\qquad j=1,\ldots,J \mbox{ and } \theta^{0}
= \bigl(\theta _{1}^{0},\ldots,\theta_{J}^{0}
\bigr).
\]
Thanks to Proposition 4.1 in~\cite{BC11}, we have:
%
\begin{prop} \label{propquad}
Suppose that Assumptions~\ref{hypP} and~\ref{hypf} hold. Then, for
any $ (\theta_{1},\ldots,\theta_{J}) \in\Theta_{\kappa} $,
\[
M(\theta_{1},\ldots,\theta_{J}) - M \bigl(
\theta_{1}^{0},\ldots,\theta _{J}^{0}
\bigr) \geq C(f,\kappa) \frac{1}{J} \sum_{j=1}^{J}
\bigl|\theta_{j} - \theta_{j}^{0}\bigr|^2,
\]
where $C(f,\kappa) = 4 \pi^{2} \vert c_1\vert^2 \cos(4 \pi\kappa) > 0$.
\end{prop}
Therefore, over the constrained set $\Theta_{\kappa}$, criterion \eqref
{defM} has a unique minimum at $\theta^{0}$ such that $M(\theta^{0}) =
0$. Let us now consider the estimators
%
\begin{equation}
\hat{\theta}^{0} =\bigl(\hat{\theta}_{1}^{0},
\ldots,\hat{\theta}_{J}^{0}\bigr) \in\argmin_{ (\theta_{1},\ldots,\theta_{J}) \in\Theta_{\kappa} }
M_{n}(\theta_{1},\ldots,\theta_{J}).
\label{eqcrittheta0}
\end{equation}
The following theorem shows that, under appropriate assumptions, the
vector $\hat{\theta}^{0}$ is a consistent estimator of $\theta^{0}$.
%
\begin{theo}\label{theoshift}
Suppose that Assumptions~\ref{hypP} and~\ref{hypf} hold. Let $J \geq
2$ and $s \geq2$. Then there exists a constant $C > 0$ that only
depends on $A$, $s$, $\kappa$, $c_*$ and $\sigma^2$ such that, for any
$f \in\tilde{W}_s(A,c_*)$, we have
%
\begin{equation}
\label{eqratenJ} \frac{1}{J} \EE^g \bigl\| \hat{
\theta}^{0}-\theta^{0}\bigr\|^2 \leq\frac
{C}{n}
\biggl(1+\frac{k_0^5}{n^{1/2}} \biggr) \biggl(1+\frac
{k_0^{3/2}J^3}{n^{1/2}} \biggr).\vadjust{\goodbreak}
\end{equation}
\end{theo}
The hypothesis $s\geq2$ in Theorem~\ref{theoshift} is related to the
need of handling the Hessian matrix associated to the criterion $M_n$.
Inequality \eqref{eqratenJ} shows that the quality of the estimation
of the random shifts depends on the ratio between $n$ and $J$. In
particular, it suggests that the quality of this estimation should
deteriorate if the number $J$ of curves increases and $n$ remains
fixed. This shows that estimating the vector $\hat{\theta}^{0} \in\RR
^{J}$ is not a standard parametric problem, since the dimension $J$ is
is allowed to grow to infinity in our setting. To the contrary, if $J$
is not too large with respect to $n$, then an estimation of the shifts\vspace*{1pt}
is feasible at the usual parametric rate $n^{-1}$. More precisely, by
Theorem~\ref{theoshift}, we immediately have the following
result:\looseness=1
%
\begin{coro} \label{coroshift}
Suppose that the assumptions of Theorem~\ref{theoshift} are satisfied.
If $k_{0} \geq1$ is a fixed integer and $J \asymp n^{\alpha}$ for some
$0 < \alpha\leq1/6$, then there exists $C_{1} > 0$ that only depends
on $A$, $s$, $\sigma^2$, $\kappa$, $c_*$ and $k_{0}$ such that
\[
\frac{1}{J} \EE^{g} \bigl\|\hat{\theta}^{0} -
\theta^{0}\bigr\|^2 \leq\frac
{C_{1} }{n}.
\]
\end{coro}
Therefore, under the additional assumption that $J \asymp n^{\alpha}$,
for some $0 < \alpha\leq1/6$, the vector $\hat{\theta}^{0}$ converges
to $\theta^{0}$ at the rate $n^{-1}$ for the normalized Euclidean norm.
The assumption $\alpha\leq1/6$ illustrates the fact that the number $J$
of curves should not be too large with respect to the size $n$ of the
design. Such a condition appears to be sufficient, but we do not claim
about the existence of an optimal rate at which the number $n$ of
design points should increase for a given increase in $J$.

\subsection{Estimation of the mean pattern}\label{subsectionestimation}
For $p \in\{0,1\}$, let $Y^{(p)}$ be given by $(Y_{q,j}^{(p)} )_{1
\leq q \leq N,  1 \leq j \leq J}$.
Thanks to the estimator $\hat{\theta}^0$ of the random shifts, we can
align the data $Y^{(1)}$ in order to estimate the mean pattern $f$ in
\eqref{eqmodelshiftedcurve}. Let $m_1<N/2$ be some positive integer.
For any $m\in\{1,\ldots,m_1\}$, we recall that the estimator $\hat
{f}_{n,J}^{(m)}$ is given by
\[
\hat{f}_{n,J}^{(m)}(t) = \frac{1}{J}\sum
_{j=1}^J\hat{f}^{(m)}_j
\bigl(t+\hat {\theta}^0_j\bigr), \qquad t\in[0,1].
\]

To simplify the notation, we omit the dependency on $k_{0}$, $n$ and
$J$ of the above estimators, and we write $\hat{f}^{(m)} = \hat
{f}_{n,J}^{(m)} $. We denote by $\EE^{(1)}$ the expectation according
to the distribution of $Y^{(1)}$. By construction, we recall that
$Y^{(0)}$ and $Y^{(1)}$ are independent. Thus, we obtain
\[
\E^{(1)} \bigl[\hat{f}^{(m)}(t) \bigr]=\bar{f}^{(m)}(t)
=\frac{1}{J}\sum_{j=1}^J
\bar{f}^{(m)}_j\bigl(t+\hat{\theta}^0_j
\bigr),\qquad  t\in[0,1],
\]
where we have set
\[
\bar{f}^{(m)}_j(t)=\sum_{\vert k\vert\leq m}
\bar{c}^{(1)}_{k,j}\phi _k(t), \qquad j\in\{1,\ldots,J\}.
\]
Therefore, $\hat{f}^{(m)}$ is a biased estimator of $f$ with respect to
$\E^{(1)}$. The idea of the procedure is that if the estimators $\hat
{\theta}^0_j$ of the shifts behave well then $d^2([f],[\bar
{f}^{(m_1)}])$ is small and estimating $f$ amounts to estimate $\bar
{f}^{(m_1)}$.

To choose an estimator of $\bar{f}^{(m_1)}$ among the $\hat
{f}^{(m)}$'s, we take a model selection approach. Before describing the
procedure, let us compute the quadratic risk of an estimator $\hat{f}^{(m)}$,
\begin{eqnarray*}
& & \E^{(1)} \biggl[\int_0^1\bigl
\llvert \bar{f}^{(m_1)}(t)-\hat {f}^{(m)}(t)\bigr\rrvert
^2\,dt \biggr]
\\
& &\qquad =\int_0^1\bigl\llvert
\bar{f}^{(m_1)}(t)-\bar {f}^{(m)}(t)\bigr\rrvert ^2\,dt+
\frac{(2m+1)\sigma^2}{NJ}.
\end{eqnarray*}
This risk is a sum of two nonnegative terms. The first one is a bias
term that is small when $m$ is close to $m_1$ while the second one is a
variance term that is small when $m$ is close to zero. The aim is to
find a trade-off between these two terms thanks to the data only. More
precisely, we choose some $\hat{m}\in\{1,\ldots,m_1\}$ such that
%
\begin{equation}
\label{msdefhatm} \hat{m}\in
\argmin_{m\in\{1,\ldots,m_1\}} \biggl\{\int_0^1\bigl
\llvert \hat {f}^{(m_1)}(t)-\hat{f}^{(m)}(t)\bigr\rrvert
^2 \,dt+\eta\frac{(2m+1)\sigma
^2}{NJ} \biggr\},
\end{equation}
where $\eta>1$ is some constant. In the sequel, the estimator that we
finally consider is $\hat{f}_{n,J}=\hat{f}^{(\hat{m})}$.

Such a procedure is well known, and we refer to Chapter 4 of \cite
{MR2319879} for more details. In particular, the estimator $\hat
{f}_{n,J}$ satisfies the following inequality:
%
\begin{eqnarray}\label{msoracle}
& & \E^{(1)} \biggl[\int_0^1\bigl
\llvert \bar{f}^{(m_1)}(t)-\hat {f}_{n,J}(t)\bigr\rrvert
^2\,dt \biggr]
\nonumber
\\
&&\qquad \leq C(\eta) \biggl\{\min_{m\in\{1,\ldots,m_1\}}\E^{(1)} \biggl[
\int_0^1\bigl\llvert \bar{f}^{(m_1)}(t)-
\hat{f}^{(m)}(t)\bigr\rrvert ^2\,dt \biggr]+\frac{\sigma^2}{NJ}
\biggr\}
\\
&&\qquad \leq C(\eta)\min_{m\in\{1,\ldots,m_1\}} \biggl\{\int_0^1
\bigl\llvert \bar {f}^{(m_1)}(t)-\bar{f}^{(m)}(t)\bigr\rrvert
^2\,dt+\frac{2(m+1)\sigma
^2}{NJ} \biggr\},\nonumber
\end{eqnarray}
where $C(\eta)>0$ only depends on $\eta$. It is known that an optimal
choice for $\eta$ is a difficult problem from a theoretical point of
view. However, in practice, taking some $\eta$ slightly greater than
$2$ leads to a procedure that behaves well as we discuss in Section \ref
{secsimus}. Moreover, for real data analysis, we often have to
estimate the variance $\sigma^2$.

\subsection{Convergence rates over Sobolev balls}

Let us denote by $\lfloor x\rfloor$ the largest integer smaller than
$x\in\RR$. We now focus on the performances of our estimation procedure
from the minimax point of view and with respect to the distance $d$
defined in \eqref{defdistance}. Note that, in Section \ref
{subsectionestimation}, we only use truncated Fourier series
expansions for building the estimators $\hat{f}^{(m)}$. In practice, we
could use other bases of $L^{2}_{\mathrm{per}}([0,1])$, and we would still have
a result like \eqref{msoracle}. In particular, the following theorem
would remain true by combining model selection techniques with bases
like piecewise polynomials or orthonormal wavelets to approximate a function.

\begin{theo} \label{theoadapt}
Assume that $nJ\geq\max\{21J,(4\sigma^2)^{2s+1}/c^{2s}\}$ where $0<c<1$
is such that $J\leq cn^{\alpha}$ for some $\alpha> 0$. Take $m_1 =
\lfloor N/2\rfloor-1$ and let $s>3/2$ and $A>0$. Then the estimator
$\hat{f}_{n,J}$ defined by procedure \eqref{msdefhatm} is such that,
for any $g\in\mathcal{G}^{\kappa}$,
\begin{eqnarray*}
& & \sup_{f\in\tilde{W}_s(A,c_*)}\mathcal{R}_g(\hat{f}_{n,J},f)
\\
& &\qquad \leq C \biggl( \vert m_1\vert^{-2s} +
m_1 n^{-2s+1} + \frac{1}{J}\EE^{g} \bigl( \bigl\|
\hat{\theta}^{0} - \theta^{0}\bigr\|^2 \bigr) + (nJ
)^{-{2s}/{(2s+1)}} \biggr)
\end{eqnarray*}
for some $C>0$ that only depends on $A$, $s$, $\sigma^2$, $\kappa$,
$k_{0}$, $\eta$ and $c$.
\end{theo}

Therefore, using the results of Corollary~\ref{coroshift} on the
convergence rate of $\hat{\theta}^0$ to~$\theta^0$, we finally obtain
the following result.

\begin{coro}\label{coroadapt}
Suppose that the assumptions of Theorems~\ref{theoshift} and~\ref{theoadapt} are satisfied. If $k_{0} \geq1$ is a fixed integer
and $J \asymp n^{\alpha}$ for some $0 < \alpha\leq1/6$, then there
exists $C' > 0$ that only depends on $A$, $s$, $\sigma^2$, $\kappa$,
$k_{0}$, $\eta$, $c_*$ and $c$ such that
\[
\sup_{g\in\mathcal{G}^{\kappa}}\sup_{f\in\tilde{W}_s(A,c_*)}\mathcal
{R}_g(\hat{f}_{n,J},f)\leq C' \bigl(
n^{-1} + (nJ )^{-
{2s}/{(2s+1)}} \bigr).
\]
\end{coro}

\section{A lower bound on the risk} \label{eqlowerbound}

The following theorem gives a lower bound on the risk over the Sobolev
ball $\tilde{W}_s(A,c_*)$.
%
\begin{theo} \label{theolowerbound}
Let us recall that
\[
\R_{n,J}\bigl(\tilde{W}_s(A,c_*),\GG^{\kappa}
\bigr) = \inf_{\hat{f}_{n,J}} \sup_{\Pr\in\GG^{\kappa}} \sup
_{f \in\tilde{W}_s(A,c_*)} \R_{\Pr}(\hat {f}_{n,J},f).
\]
There exists a constant $C > 0$ that only depends on $A$, $s$, $c_*$
and $\sigma^2$ such that
\[
\liminf_{\min(n,J) \to+ \infty} (nJ)^{{2s}/{(2s+1)}} \R_{n,J}\bigl(
\tilde {W}_s(A,c_*),\GG^{\kappa}\bigr) \geq C.
\]
\end{theo}

From the results of the previous sections, we also easily obtain the
following upperbound on the risk.

\begin{coro}
Suppose that the assumptions of Corollary~\ref{coroadapt} hold, and
assume that $2 \alpha s \leq1$.
Then, there exists a constant $C' > 0$ that only depends on $A$, $s$,
$\sigma^2$, $\kappa$, $k_{0}$, $\eta$, $c_*$ and $c$ such that
%
\begin{equation}
\sup_{\Pr\in\GG^{\kappa}} \sup_{f \in\tilde{W}_s(A,c_*)} \EE^{\Pr}
\bigl( d^{2}\bigl([\hat{f}_{n,J}],[f]\bigr) \bigr) \leq
C' (n J)^{-{2s}/{(2s+1)}}. \label{eqrateopt}
\end{equation}
\end{coro}
Note that inequality \eqref{eqrateopt} is a direct consequence of
Corollary~\ref{coroadapt} and the fact that $n^{-1} = (
{(nJ)^{-{2s}/{(2s+1)}}})$ in the settings $2 \alpha s \leq1$ and $J
\asymp n^{\alpha}$. Therefore, under the assumption that $2 \alpha s
\leq1$, the smoothed Fr\'echet mean converges at the optimal rate $(n
J)^{-{2s}/{(2s+1)}}$.

\section{Discussion} \label{secdiscussion}

As explained previously, the rate of convergence $r_{n,J} = n^{-1} + (n
J)^{-{2s}/{(2s+1)}}$ of the estimator $\hat{f}_{n,J}$ is the sum of
two terms having different interpretations. The term $(nJ)^{-
{2s}/{(2s +1)}}$ is the usual nonparametric rate that would be obtained if
the random shifts $\theta_{1}^{\ast},\ldots,\theta^{\ast}_{J}$ were
known. To interpret the second term~$n^{-1}$, let us mention the
following result that has been obtained in~\cite{BC11}.
%
\begin{prop} \label{propVanTree}
Suppose that the function $f$ is continuously differentiable. Assume
that the density $g\in\mathcal{G}^{\kappa}$ with $g(-\kappa/2)=g(\kappa
/2)=0$ and that $\mathcal{I}_g^2=\int_{-1/2}^{1/2}  ( \frac{\partial
}{\partial\theta} \log g(\theta) )^{2} g(\theta) \,d \theta<
+\infty$. Let $(\hat{\theta}_{1},\ldots,\hat{\theta}_{J})$ denote any
estimator of the true shifts $(\theta_{1}^{\ast},\ldots,\theta_{J}^{\ast
})$ computed from the $Y_{\ell,j}$'s in model \eqref
{eqmodelshiftedcurve}. Then
%
\begin{equation}
\EE^{g} \Biggl( \frac{1}{J} \sum_{j=1}^{J}
\bigl(\hat{\theta}_{j}- \theta _{j}^{\ast}
\bigr)^{2} \Biggr) \geq\frac{ \sigma^{2} }{ n \int_{0}^{1}
|f'(t)|^{2} \,dt+ \sigma^{2} \mathcal{I}_g^2}. \label{eqVanTree}
\end{equation}
\end{prop}
Proposition \eqref{propVanTree} shows that it is not possible to build
consistent estimators of the shifts by considering only the asymptotic
setting where the number of curves $J$ tends toward infinity. Indeed
inequality \eqref{eqVanTree} implies that $\liminf_{J \to+ \infty} \EE
^{g}  ( \frac{1}{J} \sum_{j=1}^{J} (\hat{\theta}_{j}- \theta
_{j}^{\ast})^{2}  ) > 0$ for any estimators $(\hat{\theta
}_{1},\ldots,\hat{\theta}_{J})$. We recall that, under the assumptions
of Corollary~\ref{coroshift}, one has
\[
\EE^{g} \Biggl( \frac{1}{J} \sum_{j=1}^{J}
\bigl|\hat{\theta}_{j}^{0} - \theta _{j}^{0}\bigr|^2
\Biggr) \leq\frac{C_{1} }{n}.
\]
The above inequality shows that, in the setting where $n$ and $J$ are
both allowed to increase, the estimation of the unknown shifts $\hat
{\theta}_{j}^{0} = \hat{\theta}_{j}^{\ast} - \frac{1}{J} \sum_{m=1}^{J}\hat{\theta}_{m}^{\ast}$
is feasible at the rate $n^{-1}$. By
Proposition~\ref{propVanTree}, this rate of convergence cannot be
improved. We thus interpret the term $n^{-1}$ appearing in the rate
$r_{n,J}$ of the smoothed Fr\'echet mean $[\hat{f}_{n,J}]$ as the price
to pay for having to estimate the shifts to compute such estimators.\vadjust{\goodbreak}

To conclude this discussion, we would like to mention the results that
have been obtained in~\cite{MR2676894} in an asymptotic setting where
only the number $J$ of curves is let going to infinity. Consider the
following model of randomly shifted curves with additive white noise:
%
\begin{eqnarray}
\label{eqmodelwhitenoise} dY_{j}(t) = f\bigl(t-
\theta_{j}^{\ast}\bigr) \,dt + \varepsilon \,dW_{j}(t),
\nonumber
\\[-8pt]
\\[-8pt]
 \eqntext{t
\in [0,1], j=1,\ldots,J \mbox{ with } \theta_{j}^{\ast}
\sim_{\mathrm{i.i.d.}} g,}
\end{eqnarray}
where the $W_{j}$'s are independent Brownian motions with $\varepsilon>
0$ being the level of additive noise. In model \eqref
{eqmodelwhitenoise}, the expectation of each observed curve $dY_{j}$
is equal to the convolution of $f$ by the density $g$ since
\[
\EE^{g} \bigl[ f\bigl(t-\theta_j^{\ast}\bigr)
\bigr] = \int f(t-\theta) g(\theta ) \,d\theta= f \star g(t).
\]
Therefore, in the ideal situation where $g$ is assumed to be known, it
has been shown in~\cite{MR2676894} that estimating $f$ in the
asymptotic setting $J \to+\infty$ (with $\varepsilon> 0$ being fixed) is
a deconvolution problem. Indeed, suppose that, for some $\nu>1/2$,
\[
\gamma_{k} = \int_{0}^{1} g(\theta)
e^{-i 2 \pi k \theta} \,d \theta \asymp|k|^{-\nu}, \qquad k \in\ZZ,
\]
with $g$ being known. Then, one can construct an estimator $\hat
{f}_{J}^{\ast}$ by a deconvolution procedure such that
\[
\sup_{f \in W_{s}(A) } \EE^{g} \int_{0}^{1}
\bigl|\hat{f}_{J}^{\ast}(t) - f(t)\bigr|^{2}\,dt \leq C
J^{-{2s}/{(2s +\nu+1)}}
\]
for some $C>0$ that only depends on $A$, $s$ and $\varepsilon$ and where
$W_{s}(A) $ is Sobolev ball of degree $s > 1/2$. Moreover, this rate of
convergence is optimal since the results in~\cite{MR2676894} show that
if $s > 2\nu+ 1$, then there exists a constant $C'>0$ that only
depends on $A$, $s$ and $\varepsilon$ such that
\[
\liminf_{J \to+ \infty} J^{{2s}/{(2s + 2 \nu+ 1)}} \inf_{\hat
{f}_{J}}
\sup_{f \in W_{s}(A) } \R(\hat{f}_{J},f) \geq C',
\]
where the above infimum is taken over the set of all estimators $\hat
{f}_{J}$ of $f$ in model~\eqref{eqmodelwhitenoise}. Hence, $r_{J} =
J^{-{2s}/{(2s + 2 \nu+ 1)}}$ is the minimax rate of convergence over
Sobolev balls in model \eqref{eqmodelwhitenoise} in the case of known
$g$. This rate is of polynomial order of the number of curves $J$, and
it deteriorates as the smoothness $\nu$ of the convolution kernel $g$
increases. This phenomenon is a well-known fact in deconvolution
problems; see, for example,~\cite{F91aos,PV99aos}. Hence, depending on
$g$ being known or not and the choice of the asymptotic setting, there
exists a significant difference in the rates of convergence that can be
achieved in a randomly shifted curves model. Our setting yields the
rate $r_{n,J} = n^{-1} + (n J)^{-{2s}/{(2s+1)}}$ (in the case where\vadjust{\goodbreak}
$J \asymp n^{\alpha}$ with $\alpha< 1/6$) that is clearly faster than
the rate $r_{J}=J^{-{2s}/{(2s + 2 \nu+ 1)}}$. Nevertheless, the
arguments in~\cite{MR2676894} also suggest that a smoothed Fr\'echet
mean in \eqref{eqmodelwhitenoise} is not a consistent estimator of $f$
if one only lets $J$ going to infinity. Therefore, the number $n$ of
design points is clearly of primary importance to obtain consistent
estimators of a mean pattern when using Fr\'echet means of curves.

\section{Numerical experiments} \label{secsimus}

\begin{figure}

\includegraphics{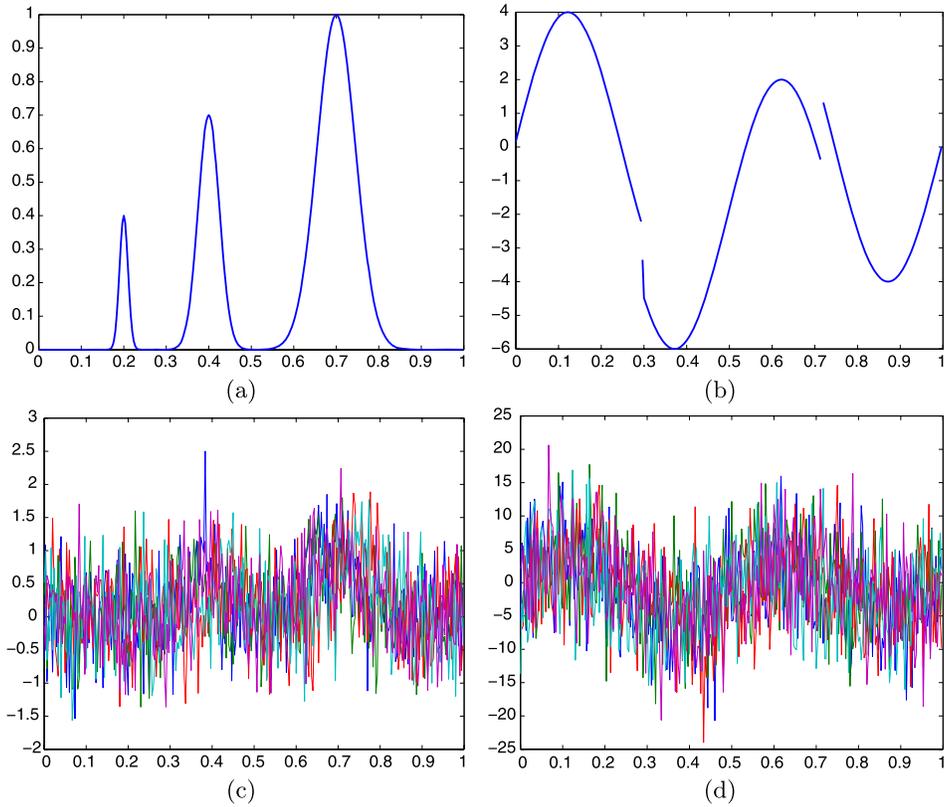}

\caption{Two test functions $f$. \textup{(a)} MixtGauss: a mixture of
three Gaussians.
\textup{(b)} HeaviSine: a piecewise smooth curve with a discontinuity.
Sample of $5$
noisy randomly shifted curves with $n=300$ for \textup{(c)} MixtGauss
and \textup{(d)} HeaviSine.} \label{Figtest}
\end{figure}

The goal of this section is to study the performance of the estimator
$\hat{f}_{n,J}$. The factors in the simulations are the number $J$ of
curves and the number $n$ of design points. As a mean pattern $f$ to
recover, we consider the two test functions displayed in Figure~\ref{Figtest}.
Then, for each combination of $n$ and $J$, we generate $M=100$
repetitions of model \eqref{eqmodelshiftedcurve} of $J$ curves with
shifts sampled from the uniform distribution on $[-\kappa,\kappa]$ with
$\kappa= 1/16$. The level of the additive Gaussian noise is measured
as the root of the signal-to-noise ratio ($\mathrm{rsnr}$) defined as
\[
\mathrm{rsnr} = \biggl( \frac{1}{\sigma^{2}}\int_{0}^{1}
\bigl(f(t) - \bar{f}\bigr)^2 \,dt \biggr)^{1/2}\qquad \mbox{where }
\bar{f} = \int_{0}^{1} f(t)\,dt,
\]
that is fixed to $\mathrm{rsnr} = 0.5$ in all the simulations. Samples of noisy
randomly shifted curves are displayed in Figure~\ref{Figtest}. For
each repetition $p\in\{1,\ldots,M\}$, we compute the estimator $\hat
{f}_{n,J,p}$ using a gradient descent algorithm to minimize criterion
\eqref{eqcrittheta0} for estimating the shifts. For all values of $n$
and $J$, we took $k_{0} = 5$ in \eqref{defMn}.
The frequency cut-off $\hat{m}$ is chosen using \eqref{msdefhatm}
with $\eta= 2.5$.

%

To analyze the numerical performance of this estimator, we have
considered the following ideal estimator that uses the knowledge of the
true random shifts $\theta_{j,p}^{\ast}$ (sampled from the $p$th replication):
\[
\tilde{f}_{n,J,p}^{(m)}(t) = \sum_{ |k| \leq m }
\Biggl( \frac{1}{J} \sum_{j=1}^{J}
\hat{c}_{k,j,p}^{(1)} e^{i 2 \pi k \theta_{j,p}^{\ast}} \Biggr)
\phi_{k}(t),\qquad t\in[0,1].
\]
The frequency cut-off $\hat{m}_{\ast}$ for the above ideal estimator is
chosen using a model selection procedure based on the knowledge of the
true shifts, that is,
\[
\hat{m}_{\ast} \in\argmin_{m\in\{1,\ldots,m_1\}} \biggl\{\int
_0^1\bigl\llvert \tilde{f}_{n,J,p}^{(m_1)}(t)-
\tilde{f}_{n,J,p}^{(m)}(t) \bigr\rrvert ^2 \,dt+\eta
\frac{(2m+1)\sigma^2}{NJ} \biggr\}
\]
with $\eta= 2.5$.

\begin{figure}

\includegraphics{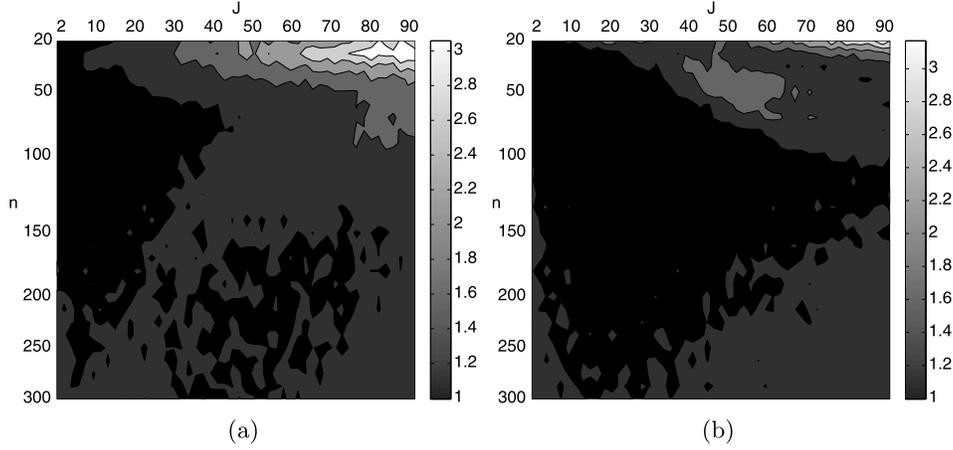}

\caption{Relative empirical error $R(n,J)$ for various values of $n$
(vertical axis) and
$J$ (horizontal axis) over $M=100$ replications: \textup{(a)} MixtGauss, \textup{(b)}
HeaviSine.} \label{figresults}
\end{figure}

Then, we define the relative empirical error between the two estimators as
\[
R(n,J) = \frac{ {{1}/{M} \sum_{p=1}^{M} d^{2}([\hat
{f}_{n,J,p}],[f])}}{  {{1}/{M} \sum_{p=1}^{M}
d^{2}([\tilde{f}_{n,J,p}^{(\hat{m}_{\ast})}],[f])}}.
\]
In Figure~\ref{figresults}, we display the ratio $R(n,J)$ for various
values of $n$ and $J$ and for the two test functions displayed in
Figure~\ref{Figtest}. It can be seen that the function $J \mapsto
R(n,J)$ is increasing. This means that the numerical performance of the
estimator $\hat{f}_{n,J}$ deteriorates as the number $J$ of curves
increases and the number $n$ remains fixed. This is clearly due to the
fact that the estimation of the shifts becomes less precise when the
dimension $J$ increases. These numerical results are thus consistent
with inequality \eqref{eqratenJ} in Theorem~\ref{theoshift} and our
discussion on the rate of convergence of $\hat{f}_{n,J}$ in Section \ref
{equpperbound}. On the other hand, the function $n \mapsto R(n,J)$ is
decreasing, and it confirms that the number $n$ of design points is a
key parameter to obtain consistent estimators of a mean pattern $f$
with Fr\'echet means of curves.

\begin{appendix}\label{app}

\section*{Appendix: Proof of the main results}

Throughout the proofs, we repeatedly use the following lemma which
follows immediately from Lemma 1.8 in~\cite{MR2724359}.\vadjust{\goodbreak}
%
\begin{lem} \label{lemmacbar}
If $f \in\tilde{W}_s(A,c_*)$, then there exists a constant $A_{0} > 0$
only depending on $A$ and $s$ such that
\[
\max_{ -{N}/{2} \leq k < {N}/{2}} \bigl|\bar{c}_{k,j}^{(p)} -
c_{k} e^{- i2\pi k \theta_{j}^{\ast}} \bigr| \leq A_{0} N^{-s + 1/2},\qquad p
\in\{ 0,1\}
\]
for all $1 \leq j \leq J$.
\end{lem}

\subsection{\texorpdfstring{Proof of Theorem \protect\ref{theoshift}}{Proof of Theorem 3.1}}

For legibility, we will write $\EE= \EE^{g}$, that is, we omit the
dependency on $g$ of the expectation. The proof is divided in several
lemmas. Let $\| \cdot\|$ denote the standard Euclidean norm in $\RR
^{J}$. First, we derive upper bounds on the second, fourth and sixth
moments of $ \| \hat{\theta}^{0}-\theta^{0}\| $. The following upper
bound on the second moment is weaker than the result that we plan to
prove. It only gives the consistency of $\hat{\theta}^0$, and we will
need some additional arguments to get the announced rate of convergence
\eqref{eqratenJ}.
%
\begin{lem}
Let $N\geq2$, $J \geq1$ and $ 1 \leq k_{0} \leq N/2$. We assume that
Assumptions~\ref{hypP} and~\ref{hypf} are satisfied, and we suppose
that $s >3/2$. Then, we have the following upper bounds, for any $f \in
\tilde{W}_{s}(A,c_*)$:
%
\begin{eqnarray}
\label{eqmoment2} \frac{1}{J} \EE \bigl(\bigl \| \hat{\theta}^{0}-
\theta^{0}\bigr\|^2 \bigr) &\leq& C_{1}
\frac{k_{0}^{1/2}}{n^{1/2}},
\\
\label{eqmoment4} \frac{1}{J} \EE \bigl( \bigl\| \hat{\theta}^{0}-
\theta^{0}\bigr\|^4 \bigr) &\leq& C_2
\frac{k_{0} J}{n}
\end{eqnarray}
and
%
\begin{equation}
\label{eqmoment6} \frac{1}{J} \EE \bigl( \bigl\| \hat{\theta}^{0}-
\theta^{0}\bigr\|^6 \bigr) \leq C_3
\frac{k_{0}^{3/2} J^2}{n^{3/2}},
\end{equation}
where $C_1, C_2$ and $C_3$ are positive constants that only depend on
$A, s, c_*, \sigma^2$ and~$\kappa$.
\end{lem}

\begin{pf}
Let $f \in\tilde{W}_s(A,c_*)$. Since $\hat{\theta}^{0} = (\hat{\theta
}_{1}^{0},\ldots,\hat{\theta}_{J}^{0})$ is a minimizer of $M_{n}$, it
follows that
\[
M\bigl(\hat{\theta}^{0}\bigr) - M\bigl(\theta^{0}\bigr)
\leq2 \sup_{\theta\in\Theta
_{\kappa}} \bigl|M_{n}(\theta) - M(\theta)\bigr|.
\]
Therefore, by Proposition~\ref{propquad}, we get
%
\begin{eqnarray}
\frac{1}{J} \EE\bigl\| \hat{\theta}^{0}-\theta^{0}
\bigr\|^2 &\leq&2 C^{-1}(c_*,\kappa) \EE \Bigl( \sup
_{\theta\in\Theta_{\kappa}} \bigl|M_{n}(\theta) - M(\theta)\bigr| \Bigr),
\label{eqbound2}
\\[-2pt]
\frac{1}{J} \EE\bigl\| \hat{\theta}^{0}-\theta^{0}
\bigr\|^4 &\leq&4 C^{-2}(c_*,\kappa) J \EE \Bigl( \sup
_{\theta\in\Theta_{\kappa}} \bigl|M_{n}(\theta) - M(\theta)\bigr|^2
\Bigr) \label{eqbound4}
\end{eqnarray}
and
%
\begin{equation}
\frac{1}{J} \EE\bigl\| \hat{\theta}^{0}-\theta^{0}
\bigr\|^6 \leq8 C^{-3}(c_*,\kappa) J^2 \EE \Bigl(
\sup_{\theta\in\Theta_{\kappa}} \bigl|M_{n}(\theta) - M(\theta)\bigr|^3
\Bigr), \label{eqbound6}
\end{equation}
where we have set $C(c_*,\kappa)=4\pi^2c_*^2\cos(8\pi\kappa)$.

Let $\theta\in\Theta_{\kappa} $ and note that $M_{n}(\theta)$ can be
decomposed as
%
\begin{equation}
M_{n}(\theta) = \bar{M}(\theta) + Q(\theta) + L(\theta),
\label{eqdecompMn}
\end{equation}
where
\begin{eqnarray*}
\bar{M}(\theta)& =& \frac{1}{J} \sum_{j=1}^{J}
\sum_{|k| \leq k_{0} } \Biggl\llvert \bar{c}_{k,j}^{(0)}
e^{2 i k \pi\theta_{j}} - \frac{1}{J} \sum_{j'=1}^{J}
\bar{c}_{k,j'}^{(0)} e^{2 i k \pi\theta_{j'}}\Biggr\rrvert
^{2},
\\[-2pt]
Q(\theta)& =& \frac{1}{N J} \sum_{j=1}^{J}
\sum_{|k| \leq k_{0} } \Biggl\llvert z_{k,j}^{(0)}
e^{2 i k \pi\theta_{j}} - \frac{1}{J} \sum_{j'=1}^{J}
z_{k,j'}^{(0)} e^{2 i k \pi\theta_{j'}}\Biggr\rrvert ^{2}
\end{eqnarray*}
and
\begin{eqnarray*}
L(\theta) & = & \frac{2}{J \sqrt{N}} \sum_{j=1}^{J}
\sum_{|k| \leq
k_{0} } \Re \Biggl[ \Biggl( \bar{c}_{k,j}^{(0)}
e^{2 i k \pi\theta_{j}} - \frac{1}{J} \sum_{j'=1}^{J}
\bar{c}_{k,j'}^{(0)} e^{2 i k \pi\theta
_{j'}} \Biggr)
\\[-2pt]
& &\hspace*{84pt}{} \times \Biggl( \overline{ z_{k,j}^{(0)}
e^{2 i k
\pi\theta_{j}} - \frac{1}{J} \sum_{j'=1}^{J}
z_{k,j'}^{(0)} e^{2 i k
\pi\theta_{j'}} } \Biggr) \Biggr].
\end{eqnarray*}

Using Lemma~\ref{lemmacbar}, it follows that, for any $\theta\in
\Theta_{\kappa}$,
\begin{eqnarray*}
& & \bigl| \bar{M}(\theta) - M(\theta) \bigr|
\\
&&\qquad\leq \frac{1}{J} \sum_{j=1}^{J}
\sum_{|k| \leq k_{0} } \Biggl\llvert \bigl| \bar {c}_{k,j}^{(0)}
\bigr| + \frac{1}{J} \sum_{j'=1}^{J} \bigl|
\bar{c}_{k,j'}^{(0)}\bigr| + 2 | c_{k} | \Biggr\rrvert
\\
& &\hspace*{61pt}\qquad{}  \times\Biggl\llvert \bigl(\bar{c}_{k,j}^{(0)}
- c_{k} e^{-2 i k
\pi\theta_{j}^{\ast}} \bigr) e^{2 i k \pi\theta_{j}}\\
&&\hspace*{77pt}\qquad{} - \frac{1}{J}
\sum_{j'=1}^{J} \bigl( \bar{c}_{k,j'}^{(0)}
- c_{k} e^{-2 i k \pi\theta
_{j'}^{\ast}}\bigr) e^{2 i k \pi\theta_{j'}}\Biggr\rrvert
\\
&&\qquad\leq 2A_{0} N^{-s + 1/2} \sum_{|k| \leq k_{0} }
\bigl( 4 | c_{k} | + 2 A_{0} N^{-s + 1/2} \bigr)
\\
&&\qquad \leq 8 A_{0}N^{-s + 1/2} (2 k_{0} + 1)
^{1/2} \sqrt{\sum_{|k| \leq
k_{0} } | c_{k}
|^{2} } + 4 A_{0}^{2} (2k_{0}+1)
N^{-2s + 1}.
\end{eqnarray*}
Hence, there exists a positive constant $C$ that only depends on $A$
and $s$ such that
%
\begin{equation}
\label{eqMdiff} \sup_{\theta\in\Theta_{\kappa}} \bigl| \bar{M}(\theta) - M(\theta) \bigr|
\leq Ck_0^{1/2}N^{-s+1/2}.
\end{equation}

Now, note that $Q(\theta) \leq\frac{ \sigma^2}{N J} Z$ with
$Z = \sum_{|k| \leq k_{0} } \sum_{j=1}^{J} \llvert z_{k,j}^{(0)} / \sigma\rrvert ^{2}$
for any \mbox{$\theta\in\Theta_{\kappa}$}. Thus, it follows that
%
\begin{equation}
\EE\sup_{\theta\in\Theta_{\kappa}}\bigl \vert Q(\theta)\bigr\vert\leq\sigma ^2
(2 k_{0} + 1) N^{-1} \label{eqQ}
\end{equation}
and
%
\begin{equation}
\EE\sup_{\theta\in\Theta_{\kappa}} \bigl\vert Q(\theta)\bigr\vert^2 \leq2
\sigma^4 (2 k_{0} + 1)^2 N^{-2}.
\label{eqQ2}
\end{equation}
By Jensen's inequality, we get
\[
\EE Z^{3/2} \leq \bigl( \EE Z^{2} \bigr)^{3/4},
\]
and since $ \EE Z^{2} \leq2 J^2 (2 k_{0} + 1)^2 $, we obtain
%
\begin{equation}
\label{eqQ32} \EE\sup_{\theta\in\Theta_{\kappa}}\bigl |Q(\theta)\bigr|^{3/2}
\leq8^{1/4} \frac{ \sigma^3}{N^{3/2}} (2 k_{0} +1)^{3/2}.
\end{equation}
Finally, using $\EE Z^{3} \leq6 J^{3} (2 k_{0} +1)^{3} $, we have
%
\begin{equation}
\label{eqQ3} \EE\sup_{\theta\in\Theta_{\kappa}} \bigl|Q(\theta)\bigr|^3 \leq6
\frac{ \sigma
^6}{N^{3}} (2 k_{0} + 1)^3.
\end{equation}

By Cauchy--Schwarz's inequality,
%
\begin{equation}
L(\theta) \leq2 \sqrt{\bar{M}(\theta)} \sqrt{Q(\theta)} \label{eqCS}.
\end{equation}
Thanks to Lemma~\ref{lemmacbar}, we get
\[
\bar{M}(\theta) \leq\frac{1}{J} \sum_{|k| \leq k_{0} }
\sum_{j=1}^{J} \bigl\llvert
\bar{c}_{k,j}^{(0)} \bigr\rrvert ^{2} \leq\sum
_{|k| \leq k_{0} } |c_{k}|^{2} +
A_{0}^2 N^{-2s + 1} (2k_{0}+1).
\]
Thus, it follows from \eqref{eqQ}, \eqref{eqQ2}, \eqref{eqQ32} and
\eqref{eqCS} that there exists a positive constant~$C'$, only
depending on $A$, $s$ and $\sigma^2$, such that
%
\begin{eqnarray}
\label{eqL} \EE\sup_{\theta\in\Theta_{\kappa}} \bigl| L(\theta) \bigr|
&\leq&
C'k_0^{1/2}N^{-1/2},
\\
\label{eqL2} \EE\sup_{\theta\in\Theta_{\kappa}} \bigl| L(\theta) \bigr|^2
&\leq&{C'}^2k_0N^{-1}
\end{eqnarray}
and
%
\begin{equation}
\label{eqL3} \EE\sup_{\theta\in\Theta_{\kappa}}\bigl | L(\theta) \bigr|^3 \leq
{C'}^3k_0^{3/2}N^{-3/2}.
\end{equation}

Since $s>3/2$, we obtain, by inequalities \eqref{eqMdiff}, \eqref
{eqQ} and \eqref{eqL},
\[
\EE \Bigl( \sup_{\theta\in\Theta_{\kappa}} \bigl|M_{n}(\theta) - M(\theta
)\bigr| \Bigr) \leq C'_1k_0^{1/2}N^{-1/2},
\]
by inequalities \eqref{eqMdiff}, \eqref{eqQ2} and \eqref{eqL2},
\[
\EE \Bigl( \sup_{\theta\in\Theta_{\kappa}} \bigl|M_{n}(\theta) - M(\theta
)\bigr|^2 \Bigr) \leq C'_2k_0N^{-1},
\]
and by inequalities \eqref{eqMdiff}, \eqref{eqQ3} and \eqref{eqL3},
\[
\EE \Bigl( \sup_{\theta\in\Theta_{\kappa}} \bigl|M_{n}(\theta) - M(\theta
)\bigr|^3 \Bigr) \leq C'_3k_0^{3/2}N^{-3/2},
\]
where $C'_1, C'_2$ and $C'_3$ are positive constants that only depend
on $A, s$ and~$\sigma^2$. Combined with inequalities \eqref
{eqbound2}, \eqref{eqbound4} and \eqref{eqbound6}, the announced
result follows from the above upper bounds.
\end{pf}

In order to prove Theorem~\ref{theoshift}, we divide the rest of the
proof in the three following steps. In the sequel of this section, we
always assume that the hypotheses of Theorem~\ref{theoshift} are
satisfied, and we use the decomposition of $M_n(\theta)$ as defined in
\eqref{eqdecompMn}.

\textit{Step} 1: there exists some positive constant $C_1$ that
only depends on $c_*$ such that
%
\begin{eqnarray}\label{propres1}
& & \frac{n}{J}\bigl\Vert\hat{\theta}^0-\theta^0
\bigr\Vert^2
\nonumber\\
&&\qquad \leq C_1 \Bigl( n J \bigl\| \grad M_{n}\bigl(
\theta^{0}\bigr) \bigr\|^2\\
&&\hspace*{48pt}{}+ nJ \sup_{\theta
\in\mathcal{U}_{\kappa}}
\bigl\llVert \grad^{2} M_{n}( \theta) -
\grad^{2} M\bigl(\theta^{0}\bigr) \bigr\rrVert
^{2}_{\mathrm{op}} \bigl\| \hat{\theta}^{0} -
\theta^{0}\bigr \|^2 \Bigr),\nonumber
\end{eqnarray}
where $\grad$ and $\grad^2$ denote the gradient and the Hessian
operators, respectively, and where we have set
\[
\mathcal{U}_{\kappa}= \bigl\{\theta\in\Theta_{\kappa}\mbox{ such
that }\bigl\Vert\theta-\theta^0\bigr\Vert\leq\bigl\Vert\hat{\theta}^0-
\theta^0\bigr\Vert \bigr\}
\]
and, for any $J\times J$ matrix $B$, the operator norm $\Vert B\Vert
_{\mathrm{op}}$ is defined by
\[
\Vert B\Vert_{\mathrm{op}}=\sup_{\theta\in\RR^J\setminus\{0\}}\frac{\Vert
B\theta\Vert}{\Vert\theta\Vert}.
\]

\textit{Step} 2: there exists some positive constant $C_2$ that
only depends on $A, s$ and $\sigma^2$ such that
%
\begin{equation}
\label{propres2} n J \EE\bigl\| \grad M_{n}\bigl(\theta^{0}
\bigr) \bigr\|^{2}\leq C_2 \biggl(1 + \frac
{k_{0}^{3}}{n}
\biggr).
\end{equation}

\textit{Step} 3: there exists some positive constant $C_3$ that
only depends on $A, s, \kappa$, $c_*$ and $\sigma^2$ such that
%
\begin{eqnarray}
\label{propres3}&& nJ \EE \Bigl(\sup_{\theta\in\mathcal{U}_{\kappa}} \bigl\llVert
\grad^{2} M_{n}( \theta) - \grad^{2} M\bigl(
\theta^{0}\bigr) \bigr\rrVert ^{2}_{\mathrm{op}} \bigl\| \hat {
\theta}^{0} - \theta^{0} \bigr\|^2 \Bigr)
\nonumber
\\[-8pt]
\\[-8pt]
\nonumber
&&\qquad\leq
C_3 \biggl(1+\frac{k_0^5}{n^{{1}/{2}}} \biggr)\frac{k_0^{{3}/{2}}J^3}{n^{{1}/{2}}}.
\end{eqnarray}

The result announced in Theorem~\ref{theoshift} follows from
inequalities \eqref{propres1}, \eqref{propres2} and~\eqref{propres3}.

\subsubsection{Proof of Step 1}

The gradients of $\bar{M}(\theta), Q(\theta)$ and $L(\theta)$ follow
from easy computations. We have, for any $1\leq\ell\leq J$,
%
\begin{eqnarray}
\label{gradmbar} \frac{\partial}{ \partial\theta_{\ell} } \bar{M}(\theta) &=& \frac{4 \pi
}{J^2} \sum
_{|k| \leq k_{0} } k \Re \Biggl[ i \overline{
\bar{c}_{k,\ell
}^{(0)} e^{2 i k \pi\theta_{\ell}}} \Biggl( \sum
_{j=1}^{J} \bar {c}_{k,j}^{(0)}
e^{2 i k \pi\theta_{j}} \Biggr) \Biggr],
\\
\label{gradq} \frac{\partial}{ \partial\theta_{\ell} } Q(\theta) &=& \frac{4 \pi
}{NJ^2} \sum
_{|k| \leq k_{0} } k \Re \Biggl[ i \overline{ z_{k,\ell
}^{(0)}
e^{2 i k \pi\theta_{\ell}} } \Biggl( \sum_{j=1,j\neq\ell}^{J}
z_{k,j}^{(0)} e^{2 i k \pi\theta_{j}} \Biggr) \Biggr]
\end{eqnarray}
and
%
\begin{eqnarray}\label{gradl}
\frac{\partial}{ \partial\theta_{\ell} } L(\theta) & = & - \frac{4 \pi
}{J^{2}\sqrt{N}} \sum
_{|k| \leq k_{0} } k \Im \Biggl[ \overline{ z_{k,\ell}^{(0)}
e^{2 i k \pi\theta_{\ell}} } \Biggl( \sum_{j=1, j\neq
\ell}^{J}
\bar{c}_{k,j}^{(0)} e^{2 i k \pi\theta_{j}} \Biggr)
\nonumber
\\[-8pt]
\\[-8pt]
\nonumber
& & \hspace*{86pt}{}- \bar{c}_{k,\ell}^{(0)} e^{2 i k \pi\theta_{\ell
}}
\Biggl( \overline{ \sum_{j=1, j\neq\ell}^{J}
z_{k,j}^{(0)} e^{2 i k
\pi\theta_{j}}} \Biggr) \Biggr].
\end{eqnarray}
Similarly, we can compute the Hessians of these functions as follows,
for $1 \leq\ell, \ell' \leq J$, if $\ell\neq\ell'$,
%
\begin{eqnarray}
\label{hessm1} \frac{\partial^2}{ \partial\theta_{\ell'}\, \partial\theta_{\ell} } \bar {M}(\theta) &=& - \frac{8 \pi^2}{J^2} \sum
_{|k| \leq k_{0} } k^2 \Re \bigl[ \overline{
\bar{c}_{k,\ell}^{(0)} } \bar{c}_{k,\ell'}^{(0)}
e^{2 i k
\pi(\theta_{\ell'}-\theta_{\ell}) } \bigr],
\\
\label{hessq1} \frac{\partial^2}{ \partial\theta_{\ell'}\, \partial\theta_{\ell} } Q(\theta)& =& -\frac{8 \pi^2}{N J^2} \sum
_{|k| \leq k_{0} } k^2 \Re \bigl[ \overline{
z_{k,\ell}^{(0)} } z_{k,\ell'}^{(0)}
e^{2 i k \pi(\theta
_{\ell'}-\theta_{\ell}) } \bigr]
\end{eqnarray}
and
%
\begin{eqnarray}\label{hessl1}
\frac{\partial^2}{ \partial\theta_{\ell'}\, \partial\theta_{\ell} } L(\theta) & = & - \frac{8 \pi^2}{J^2\sqrt{N}} \sum
_{|k| \leq k_{0} } k^2 \Re \bigl[ \bar{c}_{k,\ell'}^{(0)}
\overline{ z_{k,\ell}^{(0)} } e^{2 i k \pi(\theta_{\ell'}-\theta_{\ell}) }
\nonumber
\\[-8pt]
\\[-8pt]
\nonumber
& &\hspace*{88pt}{} + \bar{c}_{k,\ell}^{(0)} \overline{
z_{k,\ell'}^{(0)} } e^{2
i k \pi(\theta_{\ell}-\theta_{\ell'}) } \bigr],
\end{eqnarray}
and if $\ell=\ell'$,
%
\begin{eqnarray}
\label{hessm2} \frac{\partial^2}{ \partial\theta_{\ell}\, \partial\theta_{\ell} } \bar {M}(\theta) &=& \frac{8 \pi^2}{J^2} \sum
_{|k| \leq k_{0} } k^2 \Re \Biggl[ \overline{
\bar{c}_{k,\ell}^{(0)} e^{2 i k \pi\theta_{\ell}} } \Biggl( \sum
_{j=1, j\neq\ell}^{J} \bar{c}_{k,j}^{(0)}
e^{2 i k \pi\theta_{j}
} \Biggr) \Biggr],
\\
\label{hessq2} \frac{\partial^2}{ \partial\theta_{\ell}\, \partial\theta_{\ell} } Q(\theta) &=& \frac{8 \pi^2}{N J^2} \sum
_{|k| \leq k_{0} } k^2 \Re \Biggl[ \overline{z_{k,\ell}^{(0)}
e^{2 i k \pi\theta_{\ell}} } \Biggl( \sum_{j=1, j\neq\ell}^{J}
z_{k,j}^{(0)} e^{2 i k \pi\theta_{j} } \Biggr) \Biggr]
\end{eqnarray}
and
%
\begin{eqnarray}\label{hessl2}
\frac{\partial^2}{ \partial\theta_{\ell}\, \partial\theta_{\ell} } L(\theta) & = & \frac{8 \pi^2}{J^2\sqrt{N}} \sum
_{|k| \leq k_{0} } k^2 \Re \Biggl[ \bar{c}_{k,\ell}^{(0)}
e^{2 i k \pi\theta_{\ell}} \Biggl( \overline{ \sum_{j=1, j\neq\ell}^{J}
z_{k,j}^{(0)} e^{2 i k \pi\theta
_{j}} } \Biggr)
\nonumber
\\[-8pt]
\\[-8pt]
\nonumber
& &\hspace*{81pt}{}  + \overline{ z_{k,\ell}^{(0)}
e^{2 i k \pi
\theta_{\ell}} } \Biggl( \sum_{j=1, j\neq\ell}^{J}
\bar{c}_{k,j}^{(0)} e^{2 i k \pi\theta_{j}} \Biggr) \Biggr].
\end{eqnarray}

Using the fact that $\hat{\theta}^{0} \in\Theta_{\kappa}$ is a
minimizer of $M_{n}$, so $\grad M_{n}(\hat{\theta}^{0}) = 0$, a Taylor
expansion of $\theta\mapsto\grad M_n(\theta)$ with an integral form of
the remainder term leads to
%
\begin{equation}
0= \grad M_{n}\bigl(\theta^{0}\bigr) + \int
_0^1 \grad^{2} M_{n}
\bigl(\bar{\theta}(t)\bigr) \bigl(\hat{\theta}^{0} -
\theta^{0}\bigr) \,dt,
\end{equation}
where, for any $t \in[0,1]$, we have set
\[
\bar{\theta}(t)=\theta^{0}+t\bigl(\hat{\theta}^{0}-
\theta^{0}\bigr) \in\mathcal {U}_{\kappa}.
\]
Thus, we have
%
\begin{eqnarray}\label{eqhattheta}
&&\grad^{2} M\bigl(\theta^{0}\bigr) \bigl(\hat{
\theta}^{0} - \theta^{0}\bigr)
\nonumber
\\[-8pt]
\\[-8pt]
\nonumber
&&\qquad= - \grad M_{n}
\bigl(\theta^{0}\bigr) - \int_0^1
\bigl( \grad^{2} M_{n}\bigl(\bar{\theta}(t)\bigr) -
\grad^{2} M\bigl(\theta^{0}\bigr) \bigr) \bigl(\hat{
\theta}^{0} - \theta^{0}\bigr) \,dt.
\end{eqnarray}
It follows from similar computations as we did for $\bar{M}$ that
\[
\grad^{2} M\bigl(\theta^{0}\bigr) = \frac{8\pi^2}{J}
\sum_{|k| \leq k_{0}} k^2 |c_{k}|^2
\biggl( I_{J} - \frac{1}{J}\ind_J \biggr),
\]
where $I_{J}$ is the $J \times J$ identity matrix and $\ind_J$ denotes
the $J \times J$ matrix with all entries equal to one. Therefore, using
the fact that $\sum_{j = 1}^{J}  ( \hat{\theta}_{j}^{0} - \theta
_{j}^{0}  ) = 0$, we obtain
\[
\bigl\| \grad^{2} M\bigl(\theta^{0}\bigr) \bigl(\hat{
\theta}^{0} - \theta^{0}\bigr) \bigr\|^{2}=
\frac
{64\pi^4}{J^2} \biggl(\sum_{|k| \leq k_{0}} k^2
|c_{k}|^2 \biggr) ^2 \bigl\| \hat{
\theta}^{0} - \theta^{0}\bigr \|^2,
\]
and it shows that there exists a constant $C > 0$ that only depends on
$c_*$ such that
%
\begin{equation}
\bigl\| \grad^{2} M\bigl(\theta^{0}\bigr) \bigl(\hat{
\theta}^{0} - \theta^{0}\bigr) \bigr\|^{2} \geq C
\frac{1}{J^2} \bigl\| \hat{\theta}^{0} - \theta^{0}
\bigr\|^2. \label{eqhattheta2}
\end{equation}
Then, inequality \eqref{propres1} follows from \eqref{eqhattheta} and
\eqref{eqhattheta2}.

\subsubsection{Proof of Step 2}

By using Lemma~\ref{lemmacbar}, for any $1\leq k\leq k_0$ and $1\leq
\ell\leq J$, we can expand
\[
\bar{c}_{k,\ell}^{(0)} e^{2 i k \pi\theta_{\ell}^0} = c_{k}
e^{2 i k
\pi(\theta_{\ell}^0-\theta_{\ell}^\ast)} + \alpha_{k,\ell} e^{2 i k
\pi\theta_{\ell}^0}
\]
with $\vert\alpha_{k,\ell}\vert\leq A_0N^{-s+1/2}$. Because, for any
$j$, $\theta^0_j-\theta^*_j=\bar{\theta}_J$ does not depend on $j$, we have
\[
\Biggl\llvert \Re \Biggl[ i \overline{ \bar{c}_{k,\ell}^{(0)}
e^{2 i k \pi
\theta_{\ell}}} \Biggl( \sum_{j=1}^{J}
\bar{c}_{k,j}^{(0)} e^{2 i k \pi
\theta_{j}} \Biggr) \Biggr] \Biggr
\rrvert \leq A_0 J N^{-s+1/2}\bigl(2 |c_{k}| +
A_{0} N^{-s + 1/2}\bigr).
\]
Thus, by equation \eqref{gradmbar} and using Cauchy--Schwarz's
inequality, we obtain
\begin{eqnarray*}
\biggl\llvert \frac{\partial}{ \partial\theta_{\ell} } \bar{M}\bigl(\theta^{0}\bigr) \biggr
\rrvert ^{2} & \leq& \frac{32 \pi^2}{J^2} (2 k_{0} + 1)
A_0^2N^{-2s+1} \sum_{|k| \leq k_{0} }
k^2 \bigl( 4 |c_{k}|^2 + A_{0}^{2}
N^{-2s + 1} \bigr)
\\
& \leq& \frac{64 \pi^2}{J^2} (2 k_{0} + 1) A_0^2N^{-2s+1}
\bigl(A^2+A_0^2k_0^3N^{-2s+1}
\bigr).
\end{eqnarray*}
Thus, there exists a positive constant $C$ that only depends on $A$ and
$s$ such that
%
\begin{equation}
\label{eqgradM} n J\bigl \| \grad\bar{M}\bigl(\theta^{0}\bigr)
\bigr\|^{2} \leq Ck_0 \bigl(1+k_0^3n^{-2s+1}
\bigr)n^{-2s+2}.
\end{equation}

We now focus on $\grad Q$ and, by \eqref{gradq}, we can obtain
\begin{eqnarray*}
\sum_{\ell=1}^{J} \EE\biggl\llvert
\frac{\partial}{ \partial\theta_{\ell} } Q\bigl(\theta^{0}\bigr) \biggr\rrvert ^2
& \leq& \frac{16 \pi^{2}(J-1)^2}{N^2J^4} \sum_{\ell=1}^{J}
\sum_{|k| \leq k_{0} } k^2 \frac{1}{J-1} \sum
_{j=1,j\neq
\ell}^{J} \EE\bigl|z_{k,\ell}^{(0)}\bigr|^2
\bigl|z_{k,j}^{(0)}\bigr|^2
\\
& \leq& \frac{32 \pi^{2} \sigma^4k_0^3}{N^2J}.
\end{eqnarray*}
Hence, we get
%
\begin{equation}
nJ \EE\bigl\| \grad Q\bigl(\theta^{0}\bigr) \bigr\|^{2} \leq
\frac{128 \pi^{2} \sigma
^4k_0^3}{n}. \label{eqgradQ}
\end{equation}

Finally, we deal with $\grad L$. Equation \eqref{gradl} and Lemma \ref
{lemmacbar} imply
\begin{eqnarray*}
\EE\biggl\llvert \frac{\partial}{ \partial\theta_{\ell} } L(\theta) \biggr\rrvert ^{2} &
\leq& \frac{16 \pi^2 \sigma^2}{J^{4} N} \sum_{|k| \leq k_{0} } k^2
\Biggl( \Biggl\llvert \sum_{j=1, j\neq\ell}^{J}
\bar{c}_{k,j}^{(0)} e^{2 i
k \pi\theta_{j}} \Biggr\rrvert
^2 + J \bigl\llvert \bar{c}_{k,\ell}^{(0)} \bigr
\rrvert ^2 \Biggr)
\\
& \leq& \frac{16 \pi^2 \sigma^2}{J^{4} N} \sum_{|k| \leq k_{0} }
k^2 \Biggl((J-1) \sum_{j=1, j\neq\ell}^{J}
\bigl\llvert \bar{c}_{k,j}^{(0)} \bigr\rrvert ^2 +
J \bigl\llvert \bar{c}_{k,\ell}^{(0)} \bigr\rrvert ^2
\Biggr)
\\
& \leq& \frac{16 \pi^2 \sigma^2}{J^{2} N} \sum_{|k| \leq k_{0} }
k^2 \Biggl( \frac{1}{J} \sum_{j=1}^{J}
\bigl\llvert \bar{c}_{k,j}^{(0)} \bigr\rrvert ^2
\Biggr)
\\
& \leq& \frac{32 \pi^2 \sigma^2}{J^{2} N} \sum_{|k| \leq k_{0} }
k^2 \bigl( \llvert c_{k} \rrvert ^2 +
A_{0}^2 N^{-2s+1} \bigr).
\end{eqnarray*}
This last inequality leads to
%
\begin{equation}
n J \EE\bigl\| \grad L\bigl(\theta^{0}\bigr) \bigr\|^{2} \leq64
\pi^2 \sigma ^2\bigl(A^2+2A_0^2
\bigr) \bigl( 1+k_0^3 n^{-2s+1} \bigr).
\label{eqgradL}
\end{equation}
By combining inequalities \eqref{eqgradM}, \eqref{eqgradQ} and \eqref
{eqgradL}, we then obtain inequality~\eqref{propres2}.

\subsubsection{Proof of Step 3}
We introduce the Frobenius norm $\| B \|_{F}$ defined, for any $J
\times J$ matrix $B = [ B_{\ell,\ell'}]_{1 \leq\ell, \ell' \leq J}$, as
\[
\| B \|_{F} = \sqrt{\sum_{\ell, \ell' =1}^{J}
B_{\ell,\ell'}^2}.
\]
Moreover, for a self-adjoint matrix $B$, we will use the classical inequalities
%
\begin{equation}
\| B \|_{\mathrm{op}} \leq\| B \|_{F} \quad\mbox{and}\quad \| B
\|_{\mathrm{op}} \leq\max_{1 \leq\ell' \leq J} \sum
_{\ell= 1}^{J} |B_{\ell,\ell'}|.
\label{eqineqmatrix}
\end{equation}
In order to prove inequality \eqref{propres3}, we use the following
decomposition:
\begin{eqnarray}\label{eqineqHessianM}
\bigl\| \grad^{2} M_{n}( \theta) - \grad^{2} M
\bigl(\theta^{0}\bigr) \bigr\|^{2}_{\mathrm{op}} & \leq& 4
\bigl( \bigl\| \grad^{2} \bar{M}( \theta) - \grad^{2} M(\theta)
\bigr\| ^{2}_{F}
\nonumber
\\
& & \hspace*{8pt}{}+ \bigl\| \grad^{2} M( \theta) - \grad^{2} M\bigl(
\theta^{0}\bigr) \bigr\| ^{2}_{\mathrm{op}}
\\
& &\hspace*{8pt}{}+ \bigl\| \grad^{2} Q( \theta) \bigr\|^{2}_{F}+
\bigl\| \grad^{2} L( \theta ) \bigr\|^{2}_{F} \bigr).
\nonumber
\end{eqnarray}
We now deal with each term in the above inequality. Hereafter, $\ell$
and $\ell'$ always denote two integers in $\{1,\ldots,J\}$.

First, let us consider $\ell\neq\ell'$, by \eqref{hessm1} and Lemma
\ref{lemmacbar}, we get
\begin{eqnarray*}
& & \biggl\llvert \frac{\partial^2}{ \partial\theta_{\ell'}\, \partial\theta
_{\ell} } \bar{M}(\theta) - \frac{\partial^2}{ \partial\theta_{\ell'}\,
\partial\theta_{\ell} } M(
\theta) \biggr\rrvert ^2
\\
&&\qquad =  \frac{64 \pi^4}{J^4} \biggl( \sum_{|k| \leq k_{0} }
k^2 \Re \bigl[ \bigl(\bar{c}_{k,\ell}^{(0)}
\overline{ \bar{c}_{k,\ell'}^{(0)} } - \vert c_{k}
\vert^2 e^{2 i k \pi(\theta_{\ell'}^{\ast} -\theta_{\ell
}^{\ast})} \bigr)e^{i2\pi k(\theta_{\ell}-\theta_{\ell'})} \bigr]
\biggr)^2
\\
& &\qquad \leq \frac{64 \pi^4}{J^4} \biggl( \sum_{|k| \leq k_{0} }
k^2 \bigl\llvert \bar{c}_{k,\ell}^{(0)} \overline{
\bar{c}_{k,\ell'}^{(0)} } - \vert c_{k}
\vert^2 e^{2 i k \pi(\theta_{\ell'}^{\ast} -\theta_{\ell}^{\ast
})} \bigr\rrvert \biggr)^2
\\
&&\qquad \leq \frac{256 \pi^4}{J^4} \biggl( \sum_{|k| \leq k_{0} }
k^2 \bigl( A_{0} |c_{k}| N^{-s +1/2} +
A_{0}^{2} N^{-2s +1} \bigr) \biggr)^2.
\end{eqnarray*}
In the case of $\ell= \ell'$, \eqref{hessm2} and Lemma \ref
{lemmacbar} lead to
\begin{eqnarray*}
& & \biggl\llvert \frac{\partial^2}{ \partial\theta_{\ell}\, \partial\theta
_{\ell} } \bar{M}(\theta) - \frac{\partial^2}{ \partial\theta_{\ell}\,
\partial\theta_{\ell} } M(
\theta) \biggr\rrvert ^2
\\
&&\qquad \leq \frac{64 \pi^4}{J^4} \Biggl( \sum_{|k| \leq k_{0} }
k^2 \sum_{j=1, j \neq\ell}^{J} \bigl
\llvert \bar{c}_{k,\ell}^{(0)} \overline{ \bar
{c}_{k,j}^{(0)} } - \vert c_{k}\vert^2
e^{2 i k \pi(\theta_j^{\ast
}-\theta_{\ell}^{\ast}) } \bigr\rrvert \Biggr)^2
\\
&&\qquad  \leq \frac{256 \pi^4}{J^2} \biggl( \sum_{|k| \leq k_{0} }
k^2 \bigl(A_{0} |c_{k}| N^{-s +1/2} +
A_{0}^2N^{-2s +1} \bigr) \biggr)^2.
\end{eqnarray*}
Therefore, the above inequalities and \eqref{eqmoment2} imply that
there exists some positive constant $C_M$ that only depends on $A$,
$s$, $\kappa$ and $\sigma^2$ such that
%
\begin{equation}
\label{eqBoundHessM1} nJ \EE \Bigl( \sup_{\theta\in\mathcal{U}_{\kappa}} \bigl\llVert
\grad^{2} \bar{M}( \theta) - \grad^{2} M(\theta) \bigr
\rrVert ^{2}_{F} \bigl\| \hat{\theta }^{0} -
\theta^{0}\bigr \|^2 \Bigr) \leq C_M
\frac
{k_0^{6+1/2}J}{n^{2s-3/2}}.
\end{equation}

Second, using the fact that $2(1-\cos(t))\leq t^2$ for any $t\in\RR$,
if $\ell\neq\ell'$, then we have
\begin{eqnarray*}
& & \biggl\llvert \frac{\partial^2}{ \partial\theta_{\ell'}\, \partial\theta
_{\ell} } M(\theta) - \frac{\partial^2}{ \partial\theta_{\ell'}\,
\partial\theta_{\ell} } M\bigl(
\theta^{0}\bigr) \biggr\rrvert
\\
&&\qquad =  \frac{8 \pi^2}{J^2} \biggl\llvert \sum_{|k| \leq k_{0} }
k^2 |c_{k}|^2 \Re \bigl[ e^{2 i k \pi(\theta_{\ell}- \theta_{\ell}^{0} + \theta_{\ell
'}^{0} - \theta_{\ell'}) } -
1 \bigr] \biggr\rrvert
\\
& &\qquad\leq \frac{16 \pi^4}{J^2} \biggl( \sum_{|k| \leq k_{0} }
k^4 |c_{k}|^2 \biggr) \bigl|\theta_{\ell}-
\theta_{\ell}^{0} + \theta_{\ell
'}^{0} -
\theta_{\ell'}\bigr|^2,
\end{eqnarray*}
and if $\ell= \ell'$ then we obtain
\begin{eqnarray*}
& & \biggl\llvert \frac{\partial^2}{ \partial\theta_{\ell}\, \partial\theta
_{\ell} } M(\theta) - \frac{\partial^2}{ \partial\theta_{\ell}\,
\partial\theta_{\ell} } M\bigl(
\theta^{0}\bigr) \biggr\rrvert
\\[-2pt]
&&\qquad =  \frac{8 \pi^2}{J^2} \Biggl\llvert \sum_{|k| \leq k_{0} }
k^2 |c_{k}|^2 \Re \Biggl[ \Biggl( \sum
_{j=1, j\neq\ell}^{J} \bigl( e^{2 i k \pi(
\theta_{\ell}-\theta_{\ell}^{0} - \theta_{j} +\theta_{j}^{0}) } - 1
\bigr) \Biggr) \Biggr] \Biggr\rrvert
\\[-2pt]
& &\qquad \leq \frac{16 \pi^4}{J^2} \biggl( \sum_{|k| \leq k_{0} }
k^4 |c_{k}|^2 \biggr) \sum
_{j=1, j\neq\ell}^{J} \bigl|\theta_{\ell}-
\theta_{\ell
}^{0} - \theta_{j} +
\theta_{j}^{0}\bigr |^2.
\end{eqnarray*}
Therefore, by \eqref{eqineqmatrix} and under the condition $s\geq2$,
we get
\begin{eqnarray*}
\bigl\| \grad^{2} M( \theta) - \grad^{2} M\bigl(
\theta^{0}\bigr) \bigr\|_{\mathrm{op}} & \leq& \frac{32 \pi^4A^2}{J^2} \max
_{1 \leq\ell\leq J } \sum_{j=1, j\neq\ell
}^{J} \bigl|
\theta_{\ell}-\theta_{\ell}^{0} - \theta_{j}
+\theta_{j}^{0} \bigr|^2
\\[-2pt]
& \leq& \frac{64 \pi^4A^2}{J} \bigl\| \theta-\theta^{0}\bigr\|^2.
\end{eqnarray*}
Thus, by definition of $\mathcal{U}_{\kappa}$ and by \eqref
{eqmoment6}, we know that there exists some $C_M'>0$ that only depends
on $A$, $s$, $\sigma^2$ and $\kappa$ such that
%
\begin{eqnarray}\label{eqBoundHessM2}
& & nJ \EE \Bigl( \sup_{\theta\in\mathcal{U}_{\kappa}}\bigl \| \grad^{2} M(
\theta) - \grad^{2} M\bigl(\theta^{0}\bigr)
\bigr\|^{2}_{\mathrm{op}} \bigl\| \hat{\theta }^{0}-
\theta^{0}\bigr\|^2 \Bigr)
\nonumber
\\[-2pt]
&&\qquad \leq \frac{n(64\pi^4A^2)^2}{J}\EE \bigl( \bigl\| \hat{\theta}^{0}-\theta
^{0}\bigr\|^6 \bigr)
\\[-2pt]
&&\qquad  \leq C_M'\frac{k_{0}^{3/2} J^2}{n^{1/2}}.\nonumber
\end{eqnarray}

Next, we deal with the term relative to $\Vert\grad^2 Q(\theta)\Vert
_F^2$. Let us begin by noting that $ \| \hat{\theta}^{0} - \theta^{0} \|
^2 \leq4 J \kappa^2$. Thus, we have
%
\begin{eqnarray}
\label{eqboundQ1}&& \EE \Bigl( \sup_{\theta\in\mathcal{U}_{\kappa}}\bigl \|
\grad^{2} Q( \theta) \bigr\|^{2}_{F} \bigl\| \hat{
\theta}^{0} - \theta^{0} \bigr\|^2 \Bigr)
\leq 4J
\kappa^2 \EE \Bigl( \sup_{\theta\in\mathcal{U}_{\kappa}} \bigl\| \grad
^{2} Q( \theta) \bigr\|^{2}_{F} \Bigr).
\end{eqnarray}
If we take $\ell\neq\ell'$, then using \eqref{hessq1}, we get
\[
\biggl\llvert \frac{\partial^2}{ \partial\theta_{\ell'}\, \partial\theta_{\ell
} } Q(\theta) \biggr\rrvert ^2  \leq
\frac{64 \pi^{4}}{N^2 J^4} \biggl( \sum_{|k| \leq k_{0} } k^2
\bigl| z_{k,\ell}^{(0)} \bigr| \bigl| z_{k,\ell'}^{(0)} \bigr|
\biggr)^2
\]
and if $\ell=\ell'$, then by \eqref{hessq2}, we have
\begin{eqnarray*}
\llvert \frac{\partial^2}{ \partial\theta_{\ell}^2 } Q(\theta) \biggr\rrvert ^2
&\leq&
\frac{64 \pi^4}{N^2 J^4} \Biggl( \sum_{|k| \leq k_{0} } k^2
\bigl|z_{k,\ell}^{(0)}\bigr| \Biggl\llvert \sum
_{j=1, j\neq\ell}^{J} z_{k,j}^{(0)}e^{i2\pi k\theta_j}
\Biggr\rrvert \Biggr)^2
\\
&\leq& \frac{64 \pi^4}{N^2 J^3}\sum_{j=1, j\neq\ell}^J
\biggl( \sum_{|k| \leq k_{0} } k^2
\bigl|z_{k,\ell}^{(0)}\bigr| \bigl|z_{k,j}^{(0)}\bigr|
\biggr)^2.
\end{eqnarray*}
Hence, the Cauchy--Schwarz inequality leads to the following upper bound:
\begin{eqnarray*}
& & \EE \Bigl( \sup_{\theta\in\mathcal{U}_{\kappa}} \bigl\| \grad^{2} Q( \theta)
\bigr\|^{2}_{F} \Bigr)
\\
& &\qquad\leq \frac{128 \pi^4}{N^2 J^3}\EE\sum_{\ell=1}^J
\sum_{\ell'=1, \ell
'\neq\ell}^J \biggl( \sum
_{|k| \leq k_{0} } k^2 \bigl|z_{k,\ell}^{(0)}\bigr|
\bigl|z_{k,\ell'}^{(0)}\bigr| \biggr)^2
\\
& &\qquad \leq \frac{128 \pi^4}{N^2 J^3}\sum_{\ell=1}^J
\sum_{\ell'=1, \ell
'\neq\ell}^J\EE \biggl(\sum
_{\vert k\vert\leq k_0}k^2\bigl|z_{k,\ell
}^{(0)}\bigr|^2
\biggr)\EE \biggl(\sum_{\vert k\vert\leq k_0}k^2\bigl|z_{k,\ell
'}^{(0)}\bigr|^2
\biggr)
\\
&& \qquad \leq\frac{512 \pi^4\sigma^4k_0^6}{N^2 J}.
\end{eqnarray*}
Combining this bound with \eqref{eqboundQ1} gives us some $C_Q>0$ that
only depends on $\kappa$ and $\sigma^2$ such that
%
\begin{equation}
\label{eqBoundHessQ} nJ\EE \Bigl( \sup_{\theta\in\mathcal{U}_{\kappa}} \bigl\|
\grad^{2} Q( \theta) \bigr\|^{2}_{F} \bigl\| \hat{
\theta}^{0} - \theta^{0} \bigr\|^2 \Bigr) \leq
C_Q\frac{k_0^6J}{n}.
\end{equation}

Finally, we focus on the term concerning $\| \grad^{2} L( \theta) \|
^{2}_{F}$. By the Cauchy--Schwarz inequality and \eqref{eqmoment4}, we have
%
\begin{eqnarray}\label{eqLCS}
& & \EE \Bigl( \sup_{\theta\in\mathcal{U}_{\kappa}} \bigl\| \grad^{2} L( \theta)
\bigr\|^{2}_{F}\bigl \| \hat{\theta}^{0} -
\theta^{0} \bigr\|^2 \Bigr)
\nonumber
\\
&&\qquad \leq \sqrt{\EE \Bigl( \sup_{\theta\in\mathcal{U}_{\kappa}}\bigl \|
\grad^{2} L( \theta) \bigr\|^{4}_{F} \Bigr)} \sqrt{
\EE \bigl(\bigl \| \hat {\theta}^{0} - \theta^{0}
\bigr\|^4 \bigr)}
\\
& &\qquad \leq \sqrt{\EE \Bigl( \sup_{\theta\in\mathcal{U}_{\kappa}}\bigl \|
\grad^{2} L( \theta) \bigr\|^{4}_{F} \Bigr)}
\sqrt{C_2 \frac{k_{0} J^2}{n} }.\nonumber
\end{eqnarray}
Using Lemma~\ref{lemmacbar} and \eqref{hessl1}, if $\ell\neq\ell'$,
we obtain
\begin{eqnarray*}
\biggl\llvert \frac{\partial^2}{ \partial\theta_{\ell'}\, \partial\theta_{\ell
} } L(\theta) \biggr\rrvert ^2 & \leq&
\frac{64 \pi^4}{J^4 N} \biggl( \sum_{|k| \leq k_{0} } k^2
\bigl( \bigl\llvert \bar{c}_{k,\ell}^{(0)} \overline{
z_{k,\ell'}^{(0)} } \bigr\rrvert +\bigl\llvert
\bar{c}_{k,\ell'}^{(0)} \overline{ z_{k,\ell}^{(0)}
} \bigr\rrvert \bigr) \biggr)^2
\\
& \leq& \frac{64 \pi^4}{J^4 N} \biggl( \sum_{|k| \leq k_{0} }
k^2 \bigl( |c_{k} | + A_{0} N^{-s +{1}/{2}}
\bigr) \bigl(\bigl\llvert z_{k,\ell
}^{(0)}\bigr\rrvert +\bigl
\llvert z_{k,\ell'}^{(0)}\bigr\rrvert \bigr) \biggr)^2
\end{eqnarray*}
and, by \eqref{hessl2}, if $\ell= \ell'$, we get
\begin{eqnarray*}
& & \biggl\llvert \frac{\partial^2}{ \partial\theta_{\ell}\, \partial\theta
_{\ell} } L(\theta) \biggr\rrvert ^2
\\
& & \qquad\leq\frac{64 \pi^4}{J^4 N } \Biggl( \sum_{|k| \leq k_{0} }
k^2 \Biggl( \bigl|\bar{c}_{k,\ell}^{(0)}\bigr| \Biggl( \sum
_{j=1, j\neq\ell}^{J} \bigl| z_{k,j}^{(0)}
\bigr| \Biggr) + \bigl|z_{k,\ell}^{(0)}\bigr| \Biggl( \sum
_{j=1, j\neq
\ell}^{J} |\bar{c}_{k,j}^{(0)}|
\Biggr) \Biggr) \Biggr)^2
\\
& &\qquad \leq \frac{64 \pi^4}{J^4 N } \Biggl( \sum_{|k| \leq k_{0} }
k^2 \Biggl( \bigl( |c_{k} | + A_{0}
N^{-s +1/2} \bigr) \Biggl( J\bigl |z_{k,\ell
}^{(0)}\bigr| + \sum
_{j=1, j\neq\ell}^{J} \bigl| z_{k,j}^{(0)}
\bigr| \Biggr) \Biggr) \Biggr)^2.
\end{eqnarray*}
Hence, we bound the expectation in \eqref{eqLCS} from above,
\begin{eqnarray*}
& & \EE \Bigl( \sup_{\theta\in\mathcal{U}_{\kappa}} \bigl\| \grad^{2} L( \theta)
\bigr\|^{4}_{F} \Bigr)
\\
&&\qquad \leq \frac{4096\pi^8}{J^4N^2}\EE \Biggl( \Biggl[J \Biggl( \sum
_{|k| \leq
k_{0} } k^2 \bigl( |c_{k} | +
A_{0} N^{-s +1/2} \bigr) \sum_{j=1}^{J}
\bigl| z_{k,j}^{(0)}\bigr | \Biggr)^2
\\
& &\hspace*{56pt}\qquad\quad{} +\frac{4}{J}\sum_{\ell=1}^J
\biggl( \sum_{|k| \leq k_{0}
} k^2 \bigl(
|c_{k} | + A_{0} N^{-s +1/2} \bigr) \bigl\llvert
z_{k,\ell
}^{(0)}\bigr\rrvert \biggr)^2
\Biggr]^2 \Biggr)
\\
& &\qquad\leq \frac{4\times4096\pi^8}{J^2N^2}\EE \Biggl( \Biggl( \sum_{|k| \leq
k_{0} }
k^2 \bigl( |c_{k} | + A_{0} N^{-s +1/2}
\bigr) \sum_{j=1}^{J} \bigl|
z_{k,j}^{(0)}\bigr | \Biggr)^4 \Biggr)
\\
&&\qquad \leq \frac{4\times4096\pi^8}{N^2} \biggl(\sum_{|k| \leq k_{0} }
k^4 \bigl( |c_{k} | + A_{0} N^{-s +1/2}
\bigr)^2 \biggr)^2
\EE \Biggl( \Biggl( \sum_{|k| \leq k_{0} }\sum
_{j=1}^{J} \bigl| z_{k,j}^{(0)}
\bigr|^2 \Biggr)^2 \Biggr)
\\
&&\qquad \leq \frac{168\times4096\pi^8\sigma^4J^2k_0^2}{N^2} \biggl(\sum_{|k|
\leq k_{0} }
k^4 \bigl( |c_{k} | + A_{0} N^{-s +1/2}
\bigr)^2 \biggr)^2.
\end{eqnarray*}
Therefore, there exists some constant $C_L'>0$ that only depends on
$A$, $s$ and $\sigma^2$ such that
%
\begin{equation}
\label{eqLfinal} \EE \Bigl( \sup_{\theta\in\mathcal{U}_{\kappa}} \bigl\|
\grad^{2} L( \theta) \bigr\|^{4}_{F} \Bigr)\leq
C_L'\frac{k_0^2J^2}{n^2} \biggl(1 + \frac
{k_0^5}{n^{2s -1}}
\biggr)^2.
\end{equation}
Using \eqref{eqLCS} and \eqref{eqLfinal}, we know that there exists
some constant $C_L>0$ that only depends on $c_*$, $\kappa$, $A$, $s$
and $\sigma^2$ such that
%
\begin{equation}
\label{eqBoundHessL} nJ\EE \Bigl( \sup_{\theta\in\mathcal{U}_{\kappa}} \bigl\|
\grad^{2} L( \theta)\bigr \|^{2}_{F} \bigl\| \hat{
\theta}^{0} - \theta^{0} \bigr\|^2 \Bigr)\leq
C_L\frac{k_0^{3/2}J^3}{\sqrt{n}} \biggl(1 + \frac{k_0^5}{n^{2s -1}} \biggr).
\end{equation}

Finally, we use \eqref{eqBoundHessM1}, \eqref{eqBoundHessM2}, \eqref
{eqBoundHessQ} and \eqref{eqBoundHessL} with \eqref{eqineqHessianM}
to get \eqref{propres3}.

\subsection{\texorpdfstring{Proof of Theorem \protect\ref{theoadapt}}{Proof of Theorem 3.2}}

Let us assume that $f\in\tilde{W}_s(A,c_*)$. We bound the distance
between $f$ and $\hat{f}_{n,J}$ from above,
\begin{eqnarray*}
d^2\bigl([f],[\hat{f}_{n,J}]\bigr)
&=&  \inf
_{\theta\in[-1/2,1/2]}\int_0^1\bigl\llvert
f(t-\theta)-\hat{f}_{n,J}(t)\bigr\rrvert ^2\,dt
\\
& \leq& 2d^2\bigl([f],\bigl[\bar{f}^{(m_1)}\bigr]\bigr)+2\int
_0^1\bigl\llvert \bar {f}^{(m_1)}(t)-
\hat{f}_{n,J}(t)\bigr\rrvert ^2\,dt.
\end{eqnarray*}
Taking the expectation according to the distribution of $Y^{(1)}$ on
both sides and using~\eqref{msoracle} leads to
%
\begin{eqnarray}\label
{rateproofineq}
& & \E^{(1)} \bigl[d^2\bigl([f],[\hat{f}_{n,J}]
\bigr) \bigr]
\nonumber
\\
&&\qquad \leq 2d^2\bigl([f],\bigl[\bar{f}^{(m_1)}\bigr]\bigr)
\nonumber
\\
& & \qquad\quad{}+2C(\eta)\min_{m\in\{1,\ldots,m_1\}} \biggl\{\int_0^1
\bigl\llvert \bar {f}^{(m_1)}(t)-\bar{f}^{(m)}(t)\bigr\rrvert
^2\,dt+\frac{2(m+1)\sigma
^2}{NJ} \biggr\}
\\
&&\qquad \leq 2d^2\bigl([f],\bigl[\bar{f}^{(m_1)}\bigr]\bigr)
\nonumber
\\
& & \qquad\quad{}+2C(\eta)\min_{m\in\{1,\ldots,m_1\}} \Biggl\{\sum
_{m<\vert k\vert\leq
m_1}\Biggl\llvert \frac{1}{J}\sum
_{j=1}^J\bar{c}^{(1)}_{k,j}e^{i2\pi k\hat
{\theta}^0_j}
\Biggr\rrvert ^2+\frac{2(m+1)\sigma^2}{NJ} \Biggr\}.\nonumber
\end{eqnarray}

Let $\bar{\theta}_J=(\theta^*_1+\cdots+\theta^*_J)/J$, we recall that
$\theta^0=\theta^*-\bar{\theta}_J$. We begin by upper bounding the
first term. Thanks to Jensen's inequality, we obtain
%
\begin{eqnarray}\label{rateproof1}
d^2\bigl([f],\bigl[\bar{f}^{(m_1)}\bigr]\bigr) & \leq& \int
_0^1\bigl\llvert f(t-\bar{\theta
}_J)-\bar{f}^{(m_1)}(t)\bigr\rrvert ^2\,dt
\nonumber
\\
& \leq& \sum_{\vert k\vert>m_1}\vert c_k
\vert^2+\sum_{\vert k\vert\leq
m_1}\Biggl\llvert
c_ke^{-i2\pi k\bar{\theta}_J}-\frac{1}{J}\sum
_{j=1}^J\bar {c}^{(1)}_{k,j}e^{i2\pi k\hat{\theta}^0_j}
\Biggr\rrvert ^2
\\
& \leq& \sum_{\vert k\vert>m_1}\vert c_k
\vert^2+\sum_{\vert k\vert\leq
m_1}\frac{1}{J}\sum
_{j=1}^J\bigl\llvert c_ke^{-i2\pi k\bar{\theta}_J}-
\bar {c}^{(1)}_{k,j}e^{i2\pi k\hat{\theta}^0_j}\bigr\rrvert
^2.\nonumber
\end{eqnarray}
Since $f\in\tilde{W}_s(A,c_*)$, we easily upper bound the bias part
%
\begin{equation}
\label{rateproof2} \sum_{\vert k\vert>m_1}\vert
c_k\vert^2\leq A\vert m_1
\vert^{-2s}.
\end{equation}
To deal with the other part, we split it into two sums,
%
\begin{eqnarray}\label{rateproof3}
& &\frac{1}{J}\sum_{j=1}^J\bigl
\llvert c_ke^{-i2\pi k\bar{\theta}_J}-\bar {c}^{(1)}_{k,j}e^{i2\pi k\hat{\theta}^0_j}
\bigr\rrvert ^2
\nonumber
\\
&&\qquad \leq \frac{2}{J}\sum_{j=1}^J
\bigl\llvert c_ke^{-i2\pi k\bar{\theta
}_J}-c_ke^{i2\pi k(\hat{\theta}^0_j-\theta^*_j)}
\bigr\rrvert ^2
\nonumber
\\
& &\qquad\quad{} +\frac{2}{J}\sum_{j=1}^J\bigl
\llvert c_ke^{i2\pi k(\hat{\theta
}^0_j-\theta^*_j)}-\bar{c}^{(1)}_{k,j}e^{i2\pi k\hat{\theta}^0_j}
\bigr\rrvert ^2
\\
& &\qquad\leq \frac{2\vert c_k\vert^2}{J}\sum_{j=1}^J
\bigl\llvert 1-e^{i2\pi
k(\hat{\theta}^0_j-\theta^0_j)}\bigr\rrvert ^2+\frac{2}{J}\sum
_{j=1}^J\bigl\llvert \bar{c}^{(1)}_{k,j}-c_ke^{-i2\pi k\theta^*_j}
\bigr\rrvert ^2
\nonumber
\\
&&\qquad \leq \frac{8\pi^2k^2\vert c_k\vert^2}{J}\sum_{j=1}^J
\bigl(\hat {\theta}^0_j-\theta^0_j
\bigr)^2+2A_0^2N^{-2s+1},\nonumber
\end{eqnarray}
where the last inequality follows from Lemma~\ref{lemmacbar} and from
$2(1-\cos t)\leq t^2$, \mbox{$t\in\RR$}. Combining \eqref{rateproof1}, \eqref
{rateproof2} and \eqref{rateproof3}, we get, for any $g\in\mathcal
{G}^{\kappa}$,
%
\begin{eqnarray}\label{rateproofres1}
\E^g \bigl[d^2\bigl([f],\bigl[\bar{f}^{(m_1)}
\bigr]\bigr) \bigr]& \leq& A\vert m_1\vert ^{-2s}+2A_0^2(2m_1+1)N^{-2s+1}
\nonumber
\\[-8pt]
\\[-8pt]
\nonumber
& &{} +8\pi^2 \biggl(\sum_{\vert k\vert\leq m_1}k^2
\vert c_k\vert^2 \biggr)\E^g \Biggl[
\frac{1}{J}\sum_{j=1}^J \bigl(\hat{
\theta}^0_j-\theta ^0_j
\bigr)^2 \Biggr].
\end{eqnarray}

We now focus on the second term in \eqref{rateproofineq}. Let $\alpha
_{k,j} = c_ke^{-i2\pi k\theta^*_j} - \bar{c}^{(1)}_{k,j}$, using
Jensen's inequality and Lemma~\ref{lemmacbar}, for any $m\in\{1,\ldots,m_1\}$, we have
\begin{eqnarray*}
\sum_{m<\vert k\vert\leq m_1}\Biggl\llvert \frac{1}{J}\sum
_{j=1}^J\bar {c}^{(1)}_{k,j}e^{i2\pi k\hat{\theta}^0_j}
\Biggr\rrvert ^2 & \leq& \sum_{m<\vert k\vert\leq m_1}
\frac{1}{J}\sum_{j=1}^J\bigl\llvert
c_ke^{-i2\pi
k\theta^*_j}+\alpha_{k,j}\bigr\rrvert
^2
\\
& \leq& 2\sum_{\vert k\vert>m}\vert c_k
\vert^2+\frac{2}{J}\sum_{m<\vert k\vert\leq m_1}\sum
_{j=1}^J\vert\alpha_{k,j}
\vert^2
\\
& \leq& 2Am^{-2s}+4A_0^2m_1N^{-2s+1}.
\end{eqnarray*}
Let us consider $m_{\ast}$ such that
\[
m_{\ast}= \biggl\lfloor \biggl(\frac{nJ}{c} \biggr)^{1/(2s+1)}
\biggr\rfloor,
\]
where $c$ is the constant such that $J\leq cn^{\alpha}$. Note that such
a choice is allowed because it is such that $m_{\ast} \in\{1,\ldots,m_1\}$ since $n\geq21$, $s>3/2$, $\alpha\in(0,1/6]$ and $c\in(0,1)$.
In particular, such a choice leads to the following upper bound:
%
\begin{eqnarray}\label{rateproofres2}
& & \min_{m\in\{1,\ldots,m_1\}} \Biggl\{\sum_{m<\vert k\vert\leq m_1}
\Biggl\llvert \frac{1}{J}\sum_{j=1}^J
\bar{c}^{(1)}_{k,j}e^{i2\pi k\hat{\theta
}^0_j}\Biggr\rrvert
^2+\frac{2(m+1)\sigma^2}{NJ} \Biggr\}
\nonumber
\\
&& \qquad \leq4A_0^2m_1N^{-2s+1}+
\frac{2\sigma^2}{NJ}+2\min_{m\in\{1,\ldots,m_1\}} \biggl\{Am^{-2s}+
\frac{m\sigma^2}{NJ} \biggr\}
\nonumber
\\
&&\qquad \leq 4A_0^2m_1N^{-2s+1}+
\frac{2\sigma^2}{NJ}
\\
& & \qquad\quad{}+2 \biggl(\frac{A}{2}c^{2s/(2s+1)}+2\sigma^2c^{2s/(2s+1)}
\biggr) (nJ )^{-2s/(2s+1)}
\nonumber
\\
&&\qquad \leq 4A_0^2m_1N^{-2s+1}+ \bigl(1+
\bigl(A+4\sigma^2\bigr)c^{2s/(2s+1)} \bigr) (nJ )^{-2s/(2s+1)}.\nonumber
\end{eqnarray}
Putting \eqref{rateproofres1} and \eqref{rateproofres2} in \eqref
{rateproofineq} leads to, for any $g\in\mathcal{G}^{\kappa}$,
\begin{eqnarray*}
\R_g(\hat{f}_{n,J},f) & \leq& A\vert m_1
\vert^{-2s}+4A_0^2\bigl(\bigl(1+C(\eta )
\bigr)m_1+1\bigr)N^{-2s+1}
\\
& &{} +8\pi^2 \biggl(\sum_{\vert k\vert\leq m_1}k^2
\vert c_k\vert^2 \biggr)\E^g \Biggl[
\frac{1}{J}\sum_{j=1}^J \bigl(\hat{
\theta}^0_j-\theta ^0_j
\bigr)^2 \Biggr]
\\
& & {} + C(\eta) \bigl(1+\bigl(A+4\sigma^2\bigr)c^{2s/(2s+1)}
\bigr) (nJ )^{-2s/(2s+1)},
\end{eqnarray*}
that completes the proof using the fact that $m_1 = \lfloor N/2\rfloor-1$.

\subsection{\texorpdfstring{Proof of Theorem \protect\ref{theolowerbound}}{Proof of Theorem 4.1}}

The arguments that we use to derive this result are based on Assouad's
cube lemma; see, for example,~\cite{MR2724359}. This lemma is
classically used in nonparametric statistics to derive lower bounds on
a risk. We will show that one can construct a set of functions $\FF_{0}
\subset\tilde{W}_s(A,c_*)$ such that there exists a constant $C >0$
(only depending on $A$, $s$, $c_*$ and $\sigma^2$) such that, for any
large enough $n$ and $J$,
%
\begin{equation}
\qquad\R_{n,J}\bigl(\tilde{W}_s(A,c_*),\GG^{\kappa}
\bigr) \geq\inf_{\hat{f}_{n,J}} \sup_{\Pr\in\GG^{\kappa}} \sup
_{f \in\FF_{0}} \R_g(\hat{f}_{n,J},f) \geq C
(nJ)^{-{2s}/{(2s+1)}}, \label{eqminRisk0}
\end{equation}
where $\hat{f}_{n,J}$ denote some estimator of $f$. For the sake of
legibility, we assume in the sequel that $c_*=1$. Let
\[
\FF_{0} = \biggl\{f_{w} \dvtx t\in[0,1] \mapsto\sqrt{
\mu_{n,J}} \sum_{k \in
K_{n,J} } w_{k}
\phi_{k}(t),w_{k} \in\{-1,1\},w_{-k} =
w_{k} \biggr\},
\]
where $K_{n,J} =  \{k \in\ZZ,   0 < |k| \leq D_{n,J}  \}$,
$\mu_{n,J}$ is a positive real and $D_{n,J}$ is a positive integer that
will be specified below. Let us introduce the notation $\Omega= \{
-1,1\}^{D_{n,J}}$ and note that any function $f_{w}\in\mathcal{F}_0$ is
parametrized by a unique element $w \in\Omega$. Under the condition
%
\begin{equation}
\mu_{n,J} = c D_{n,J}^{-2s-1}\qquad \mbox{with } c \leq
A, \label{eqmu}
\end{equation}
it can easily be checked that $\FF_{0} \subset\tilde{W}_s(A,c_*)$. In
what follows, $D_{n,J}$ is chosen as the largest integer smaller that
$(nJ)^{{1}/{(2s+1)}}$. Hereafter, $\EE_{w}^{g}$ will denote the
expectation with respect to the distribution $\P_{w}^{g}$ of the random
vector $(Y_{\ell,j})_{1 \leq\ell\leq n, 1 \leq j \leq J} \in\RR
^{nJ}$ in model \eqref{eqmodelshiftedcurve} under the hypothesis that
$f = f_{w}$ and the assumption that the shifts are i.i.d. random
variables with density $g \in\GG^{\kappa}$. Note that for any $g \in
\GG^{\kappa}$
\begin{eqnarray*}
\sup_{f \in\FF_{0}} \R_{g}(\hat{f}_{n,J},f) & = &
\sup_{f \in\FF_{0}} \EE \biggl[ \inf_{ \theta\in[0,1] } \biggl(
\int_{0}^{1} \bigl|\hat {f}_{n,J}(t-\theta) -
f(t)\bigr|^{2}\,dt \biggr) \biggr]
\\
& \geq& \frac{1}{|\Omega|} \sum_{w \in\Omega}
\EE^{g}_{w} \biggl[ \inf_{ \theta\in[0,1] } \biggl(
\int_{0}^{1} \bigl|\hat{f}_{n,J}(t-\theta) -
f_{w}(t)\bigr|^{2}\,dt \biggr) \biggr]
\\
& \geq& \frac{1}{|\Omega|} \sum_{w \in\Omega}
\EE^{g}_{w} \biggl[ \inf_{ \theta\in[0,1] } \sum
_{k \in K_{n,J} } \bigl\llvert \hat{c}_{k}
e^{-2 i k
\pi\theta}- \sqrt{\mu_{n,J}} w_{k} \bigr\rrvert
^2 \biggr],
\end{eqnarray*}
where $\hat{c}_{k} = \int_{0}^{1} \hat{f}_{n,J}(t) \overline{\phi
_{k}(t)} \,dt $ is the $k$th Fourier coefficient of $\hat{f}_{n,J}$. Now,
we consider, for $k \in K_{n,J} $ and $\theta\in[0,1]$,
\[
\hat{w}_{k, \theta} \in\argmin_{v \in\{-1,1\}} \bigl|\hat{c}_{k}
e^{-2 i k
\pi\theta} - \sqrt{\mu_{n,J}}v\bigr|.
\]
We have the inequality
\begin{eqnarray*}
\bigl|\sqrt{\mu_{n,J}}\hat{w}_{k, \theta} - \sqrt{\mu_{n,J}}
w_{k} \bigr| & \leq& \bigl|\hat{c}_{k} e^{-2 i k \pi\theta} - \sqrt{
\mu_{n,J}}\hat {w}_{k, \theta} \bigr|
\\
& & {}+ \bigl|\hat{c}_{k} e^{-2 i k \pi\theta} - \sqrt{\mu_{n,J}}
w_{k} \bigr|
\\
& \leq& 2 \bigl|\hat{c}_{k}e^{-2 i k \pi\theta} - \sqrt{\mu_{n,J}}
w_{k}\bigr |
\end{eqnarray*}
that implies
\begin{eqnarray*}
\sup_{f \in\FF_{0}} \R_{g}(\hat{f}_{n,J},f) &
\geq& \frac{\mu
_{n,J}}{4 |\Omega| } \sum_{k \in K_{n,J} } \sum
_{w \in\Omega} \EE ^{g}_{w} \Bigl[ \inf
_{ \theta\in[0,1] } \bigl( \llvert \hat{w}_{k,
\theta} -
w_{k} \rrvert ^2 \bigr) \Bigr].
\end{eqnarray*}
For $w \in\Omega$ and $k \in K_{n,J}$, we define $w^{(k)} \in\Omega$
such that, for any $\ell\neq k$, $w^{(k)}_{\ell} = w_{\ell}$ and
$w^{(k)}_{k} = -w_{k}$. Then it follows that
%
\begin{equation}
\sup_{f \in\FF_{0}} \R_{g}(\hat{f}_{n,J},f)
\geq \frac{\mu
_{n,J}}{4 |\Omega| } \sum_{k \in K_{n,J} } \sum
_{w \in\Omega| w_{k} =
1} R_{k},\label{eqmin00}
\end{equation}
where we have set
\[
R_{k} = \EE^{g}_{w} \Bigl[ \inf
_{ \theta\in[0,1] } \bigl( \llvert \hat {w}_{k, \theta} -
w_{k} \rrvert ^2 \bigr) \Bigr] + \EE^{g}_{w^{(k)}}
\Bigl[ \inf_{ \theta\in[0,1] } \bigl( \llvert \hat{w}_{k, \theta} +
w_{k} \rrvert ^2 \bigr) \Bigr].
\]
Let $\theta^{\ast} = (\theta_{1}^{\ast},\ldots,\theta_{J}^{\ast})$, we
introduce the notation $\EE_{w}^{\theta^{\ast}}$ to denote expectation
with respect to the distribution $\P_{w}^{\theta^{\ast}}$ of the random
vector $(Y_{\ell,j})_{1 \leq\ell\leq n, 1 \leq j \leq J} \in\RR
^{nJ}$ in model \eqref{eqmodel} conditionally to $\theta_{1}^{\ast
},\ldots,\theta_{J}^{\ast}$.
Hence, using this notation, we have
%
\begin{equation}
\label{eqRk} R_{k} = \int_{[-{1}/{2},{1}/{2}]^{J} }
R_{k}\bigl(\theta^{\ast}\bigr) g\bigl(\theta_{1}^{\ast}
\bigr) \cdots g\bigl(\theta_{J}^{\ast}\bigr) \,d
\theta_{1}^{\ast
}\cdots \,d\theta_{J}^{\ast},
\end{equation}
where
\[
R_{k}\bigl(\theta^{\ast}\bigr) = \EE^{\theta^{\ast}}_{w}
\Bigl[ \inf_{ \theta\in
[0,1] } \bigl( \llvert \hat{w}_{k, \theta} -
w_{k} \rrvert ^2 \bigr) \Bigr] + \EE^{\theta^{\ast}}_{w^{(k)}}
\Bigl[ \inf_{ \theta\in[0,1]
} \bigl( \llvert \hat{w}_{k, \theta} +
w_{k} \rrvert ^2 \bigr) \Bigr].
\]
Now, note that for any $0 < \delta< 1$,
%
\begin{eqnarray}\label{eqmin01}
R_{k}\bigl(\theta^{\ast}\bigr) & = & \EE^{\theta^{\ast}}_{w}
\biggl[ \inf_{ \theta
\in[0,1] } \bigl( \llvert \hat{w}_{k,\theta} -
w_{k} \rrvert ^2 \bigr) + \inf_{ \theta\in[0,1] }
\bigl( \llvert \hat{w}_{k, \theta} + w_{k} \rrvert ^2
\bigr) \frac{d \P^{\theta^{\ast}}_{w^{(k)}}}{d \P^{\theta
^{\ast}}_{w}}(Y) \biggr]
\nonumber
\\[-8pt]
\\[-8pt]
\nonumber
& \geq& 4 \EE^{g}_{w} \min \biggl(1,
\frac{d \P^{\theta^{\ast
}}_{w^{(k)}}}{d \P^{\theta^{\ast}}_{w}}(Y) \biggr) \geq4 \delta\P ^{\theta^{\ast}}_{w}
\biggl( \frac{d \P^{\theta^{\ast}}_{w^{(k)}}}{d \P
^{\theta^{\ast}}_{w}}(Y) \geq\delta \biggr),
\end{eqnarray}
where $Y \in\RR^{nJ}$ is the random vector obtained from the
concatenation of the observations from model \eqref
{eqmodelshiftedcurve} under the hypothesis $f = f_{w}$ and
conditionally to $\theta_{1}^{\ast},\ldots,\theta_{J}^{\ast}$. Because
$ w_{k} = 1$, we know that
\[
\log\frac{d \P^{\theta^{\ast}}_{w^{(k)}}}{d \P^{\theta^{\ast
}}_{w}}(Y)=- \frac{2}{\sigma^2} \sum
_{j= 1}^{J} \sum_{\ell= 1}^{n}
\mu _{n,J} \bigl|\phi_{k}\bigl(t_{\ell}-
\theta_{j}^{\ast}\bigr)\bigr|^2 + \sqrt{\mu
_{n,J}}\varepsilon_{\ell,j} \phi_{k}
\bigl(t_{\ell}-\theta_{j}^{\ast}\bigr).
\]
Therefore, $\log\frac{d \P^{\theta^{\ast}}_{w^{(k)}}}{d \P^{\theta
^{\ast}}_{w}}(Y) $ is a random variable that is normally distributed
with mean $- \frac{2}{ \sigma^2} nJ \mu_{n,J}$ and variance $ \frac{4}{
\sigma^2} nJ \mu_{n,J}$. Now, since $D_{nJ}$ is the largest integer
smaller than $(nJ)^{{1}/{(2s +1)}}$, it follows from equation \eqref
{eqmu} that, for any $n$ and $J$ large enough,
\[
0 \leq nJ \mu_{n,J} \leq2 A.
\]
Thus, there exists $0 < \delta< 1$ and a constant $c_{\delta} > 0$
(only depending on $A$, $\sigma^2$ and $\delta$) such that
\[
\P^{\theta^{\ast}}_{w} \biggl( \frac{d \P^{\theta^{\ast}}_{w^{(k)}}}{d \P
^{\theta^{\ast}}_{w}}(Y) \geq\delta
\biggr) \geq c_{\delta}.
\]
Combining this inequality with \eqref{eqmin00}, \eqref{eqRk} and
\eqref{eqmin01} leads to
%
\begin{equation}
\sup_{f \in\FF_{0}} \R_{g}(\hat{f}_{n,J},f) \geq
\frac{4\delta\mu
_{n,J}}{|\Omega| } \sum_{k \in K_{n,J} } \sum
_{w \in\Omega| w_{k} =
1} c_{\delta} \geq\delta c_{\delta}
\mu_{n,J} D_{n,J}. \label{eqmin03}
\end{equation}
Since $\mu_{n,J} = c D_{n,J}^{-2s-1}$ and $D_{n,J} \leq(nJ)^{{1}/{(2s+1)}}$, it follows that
\[
\mu_{n,J} D_{n,J} = c D_{n,J}^{-2s} \geq
c (nJ)^{-{2s}/{(2s+1)}},
\]
which combined with \eqref{eqmin03} proves inequality \eqref
{eqminRisk0} and completes the proof of Theorem~\ref{theolowerbound}.
\end{appendix}

\section*{Acknowledgments}
We are very much indebted to the referees and the Associate Editor
for their constructive comments and remarks that helped us
to improve the presentation of the original manuscript.


%

\printaddresses

\end{document}